%% file: main_arxiv.tex
\documentclass[final,onefignum,onetabnum]{siamart190516}

\usepackage{tikz,pgfplots}
\pgfplotsset{width=10cm,compat=1.9,
}

\input{siopt_shared}

\begin{document}

\maketitle

\begin{abstract}
We consider the problem of decomposing a regular non-negative function as a sum of squares of functions which preserve some form of regularity. In the same way as decomposing non-negative polynomials as sum of squares of polynomials allows to derive methods in order to solve global optimization problems on polynomials, decomposing a regular function as a sum of squares allows to derive methods to solve global optimization problems on more general functions. As the regularity of the functions in the sum of squares decomposition is a key indicator in analysing the convergence and speed of convergence of optimization methods, it is important to have theoretical results guaranteeing such a regularity. 
In this work, we show second order sufficient conditions in order for a $p$ times continuously differentiable non-negative function to be a sum of squares of $p-2$ differentiable functions. The main hypothesis is that, locally, the function grows quadratically in directions which are orthogonal to its set of zeros. The novelty of this result, compared to previous works is that it allows sets of zeros which are continuous as opposed to discrete, and also applies to manifolds as opposed to open sets of $\R^d$. This has applications in problems where manifolds of minimizers or zeros typically appear, such as in optimal transport, and for minimizing functions defined on manifolds.
\end{abstract}



\newcommand{\noteuly}[1]{{\bfseries \color{red} #1}}

\newcommand{\notefb}[1]{{\bfseries \color{red} #1}}

\section{Introduction}

The relationship between non-negative functions and  functions decomposable as sums of squares is a fundamental question in both theoretical and applied mathematics. 
From a theoretical viewpoint, the decomposability of a non-negative function in terms of sum of squares is the basis of important theoretical objects and properties: quadratic modules \cite{marshall2006} in algebraic geometry, regularizing operators such as Laplacians or sub-Laplacians in (sub-)Riemannian geometry \cite{Hormander1967,bony1998}, non-negative symbols in pseudo-differential calculus \cite{Hormander2007,tataru2002}.
From an applicative viewpoint, representing a non-negative function in terms of sum of squares allows to simplify the analysis of probability representations and optimization problems \cite{lasserre2010moments,marteau2020non}. Restricting to the case of  non-negative polynomials this has been applied to global optimization and generalized methods of moments \cite{lasserre2010moments,Lasserre2020}. More generally, the decomposition of non-negative $p$-times differentiable functions allowed to derive simple and fast optimization algorithms in the context of global optimization \cite{rudi2020finding}, the Kantorovich problem in optimal transport \cite{vacher21a}, some formulations of optimal control \cite{berthier2021}. Moreover, it allowed to obtain an effective and concise representation for probability densities, with applications in probabilistic inference, sampling, machine learning \cite{rudi2021psd,marteauferey2021sampling}.

\subsection*{The importance of preserving regularity}

In this work, we state sufficient conditions for a non-negative function $f$ to be written as a sum of squares of functions $f_i$. Of course, if no other constraints are added, this is a trivial problem as writing $f = (\sqrt{f})^2$ would offer an immediate solution. What we want to understand here are sufficient conditions which allow to inherit a form of regularity of the function $f$ in the sum of squares decomposition. This is not necessarily the case when taking the square root: for example, the map $(x,y) \mapsto x^2 + y^2$ is smooth and a sum of smooth squares, but its square root is not differentiable at $(0,0)$.
 
  In general, being able to decompose the function with a certain regularity is important. Of course, there is a complex interaction between the structural constraint of being a sum of squares and the original regularity of the function, and the two may not work very well together (see \cref{thm:impossibility}). However, for certain theoretical and applied problems, it is crucial to maintain some regularity. For example, in the setting introduced in \cite{rudi2020finding}, the speed of convergence of the presented algorithm of global optimization depends on the regularity of the sum of squares representation of the function $f - f_*$ (where $f_*$ is the global minimum of $f$). In polynomial sum-of-squares (SoS) optimization, the running time depends on the degree in the sum of squares decomposition of $P - P_*$.

 Abstractly, we can formulate the following generic question. If $f \in \cc_1$ where $\cc_1$ describes a form of regularity, can we write $f$ as a sum of squares of functions of class $\cc_2$, where $\cc_2$ inherits the regularity properties $\cc_1$ as much as possible?
 
\subsection*{Problem setting}







In this work, we will concentrate on the class  $C^p$ of $p$ times differentiable functions with continuous $p$-th derivatives on $\R^d$ (or any $d$-dimensional  manifolds $M$). For simplicity, in this introduction, we will state the main results for functions on $\R^d$. We will show that under a certain condition on the set of zeros $\Z$ of $f$, if $f$ is a $C^p$ non-negative function on $\R^d$, it can be decomposed as
\begin{equation}
    \label{eq:sos_form}
    f = \sum_{i  \in I}{f_i^2},\qquad f_i \in C^{p-2}(\R^d),
\end{equation}
 where $(f_i)$ is an at most countable family and has locally finite support.
Two elements are important in \cref{eq:sos_form}: the locally finite aspect and the regularity of the functions $f_i$, i.e., $p-2$. This is a consequence of the fact that we will consider \emph{second order sufficient conditions}, hence the loss of two derivatives.  

\subsection{Intuition and previous results} Let us give an intuition as to how we obtain decompositions in the form \cref{eq:sos_form}. First, we start by proving that this decomposition holds locally in a neighborhood of any $x_0 \in \R^d$. It is then possible to invoke a result to ``glue'' the local decompositions together; we develop the tools to do so in \cref{sec:gluing} (note that this is one of the key differences between results for polynomials and results for functions). For any fixed $x_0 \in \R^d$, if $f(x_0) >0$, then $f_1:= \sqrt{f}$ is well defined and of class $C^p$ around $x_0$, and so \cref{eq:sos_form} holds locally around $x_0$ since $f = f_1^2$. The crux of the problem is to determine whether $f$ can be decomposed as a sum of squares around a point in the set of zeros $\Z$ of $f$, i.e., the set of points $x$ such that $f(x)=0$. Since $f$ is non-negative, all such points are necessarily minimizers of $f$, hence the following \emph{necessary second-order condition}:
\begin{equation}
    \label{eq:necessary_condition}
    \forall x_0 \in \Z,~ \nabla f (x_0)= 0,~\nabla^2 f(x_0) \succeq 0.
\end{equation}

Around any $x_0 \in \Z$, $f$ can be approximated by a parabola since the eigenvalues of $\nabla^2 f(x_0)$ are non-negative: $f(x) = x^\top \nabla^2 f(x_0) x + o(\|x\|^2)$ using a Taylor expansion. Since any parabola can be written as the sum of at most $d$ squares of linear functions (just write the eigen-decomposition of $\nabla^2 f(x_0)$), we see that up to the $o(\|x\|^2)$ factor, we can indeed write $f$ as a sum of at most $d$ squares around $x_0$. The whole difficulty of the following results is to go beyond this $o(\|x\|^2)$ approximation and have an exact decomposition, using the Taylor expansion with integral remainder.

It turns out that in the case where $\nabla^2 f(x_0) \succ 0$, that is when the Hessian has strictly positive eigenvalues, this decomposition can be made exact. We will call this condition the {\em strict Hessian condition} (SHC) at $x_0$. This result exists in recent work: it is a particular case of Theorem 2 of \cite{rudi2020finding}, applied to the set ${\cal H} = C^{p-2}$. Precisely, it states

\begin{theorem}[Theorem 2 of \cite{rudi2020finding}]\label{thm:simple_minima}
Let $f$ be a non-negative function of class $C^p$ for $p \geq 2$, and assume that the zeros $\Z$ of $f$ satisfy the strict Hessian condition:
\begin{equation}\label{eq:sc_1}
    \forall x_0 \in \Z,~ \nabla^2 f(x_0) \succ 0.
\end{equation}
If $f$ has a finite number $m = |\Z|$ of zeros, then $f$ satisfies \cref{eq:sos_form} with $dm + 1$ functions $f_i$.
\end{theorem}

This situation is illustrated on the left hand side of \cref{fig:example_intro}, where the Hessian is positive definite at all four zeros of $f$ and hence satisfies the SHC: by \cref{thm:simple_minima}, it can be decomposed as a sum of squares. It is not the case on the right hand side, where there is a continuous subspace of zeros: in that case, $f$ does not satisfy the SHC. 

\subsection*{Contribution}

\begin{figure}\label{fig:example_intro}
\begin{tikzpicture}[scale = 0.7]
\begin{axis}
  \node at (axis cs:1,1,0) [circle, scale=0.7, draw=black!100,fill=black!100] {};
  \node at (axis cs:1,-1,0) [circle, scale=0.7, draw=black!100,fill=black!100] {};
  \node at (axis cs:-1,1,0) [circle, scale=0.7, draw=black!100,fill=black!100] {};
  \node at (axis cs:-1,-1,0) [circle, scale=0.7, draw=black!100,fill=black!100] {};
\addplot3[
    opacity=0.4,
    surf,
    domain=-2.2:2.2,
    y domain=-2.2:2.2,
]
{((x-1)^2 + (y-1)^2)*((x-1)^2 + (y+1)^2)*((x+1)^2 + (y-1)^2)*((x+1)^2 + (y+1)^2)*exp(-4 *ln(3+x^2 + y^2))};
\end{axis}
\end{tikzpicture}
\begin{tikzpicture}[scale = 0.7]
\begin{axis}
\addplot3[
    surf,
    opacity = 0.4,
    domain=-pi:pi,
    y domain=-pi:pi
]
{((x - 1- cos(180*y/pi))^2)*((x+2)^2 + y^2)*exp(-ln(1+(x+2)^2 + y^2))};
\addplot3[
    line width = 3pt,
    domain=-pi:pi,
    samples = 100,
    samples y=0,
]
({1+cos(180*x/pi)},
{x},
{0});
\node at (axis cs:-2,0,0) [circle, scale=0.7, draw=black!100,fill=black!100] {};
\end{axis}
\end{tikzpicture}

\vspace*{-.25cm}

\caption{Plots of functions $z = f(x,y)$, where the zeros of $f$ are highlighted in black. \textbf{left}: $f$ satisfies the SHC, \textbf{right}: $f$ satisfies the NHC but not the SHC.}
\end{figure}
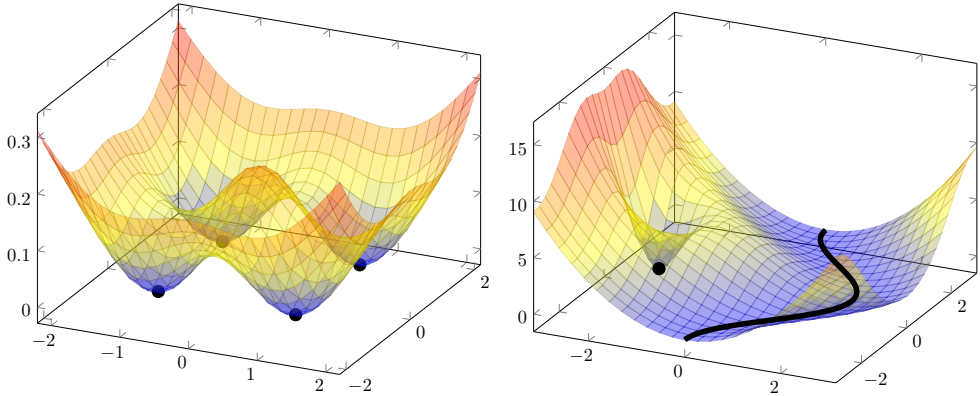

While the SHC  condition \cref{eq:sc_1} already offers a nice result in \cref{thm:simple_minima}, we see that there is a big difference with the necessary condition \cref{eq:necessary_condition}. Previous results in the literature show that \cref{eq:necessary_condition} is not sufficient to be decomposed as a sum of squares of $C^{p-2}$ functions as soon as the dimension $d$ is greater than $3$ 
(see \cref{thm:impossibility} in the background section form more details). On the other hand, \cref{eq:sc_1} is very restrictive. In particular, it implies that the set $\Z$ of zeros is discrete. However, in some situations such that of \cite{vacher21a}, the set of zeros has a natural structure, which can be a sub-manifold of $\R^d$ (consider for instance the extreme case where $f=0$). In this paper, we show that if the set $\Z$ of zeros is a sub-manifold of $\R^d$ such that the Hessian of $f$ along this manifold is positive along all directions which are not tangent to $\Z$, then \cref{eq:sos_form} still holds. This is the case for the function depicted in the right hand side of \cref{fig:example_intro}, and illustrates the difference between previous works and our contributions. More formally, we prove the following result.

\begin{theorem}\label{thm:general_Rd_intro}
Let $f$ be a non-negative function of class $C^p$ for $p \geq 2$ and let $\Z$ denote the set of zeros of $f$. If $\Z$ is a sub-manifold of $\R^d$ of class $C^1$ such that 
\begin{equation}\label{eq:first_nhc}
    \forall x_0 \in \Z,~ \forall h \in \R^d \setminus{T_{x_0}\Z},~ h^{\top} \nabla^2 f(x_0) h > 0,
\end{equation}
then $f$ satisfies \cref{eq:sos_form}, and $\Z$ is of class $C^{p-1}$. Here, $T_{x_0} \Z$ denotes the tangent space to $\Z$ at $x_0$, which is a vector sub-space of $\R^d$.
\end{theorem}

This theorem is proved as \cref{thm:main_Rd} in \cref{sec:euclidean}, and the assumption \cref{eq:first_nhc} will be referred to as the normal Hessian condition (or NHC). Note that the NHC assumption encompasses that of the SHC assumption of \cref{thm:simple_minima}; in that case, the results presented in this paper make the result tighter by removing the assumption that $\Z$ be finite and by needing only $d + 1$ squares to represent the function, and not $d|\Z| + 1$ (see the full version of \cref{thm:main_Rd}).

The proof techniques used to prove this theorem differ from the proof of \cite{rudi2020finding} and use tools from differential geometry and Morse theory. In particular, the proof extends naturally to functions defined on $d$-dimensional manifolds, which is the object of \cref{sec:manifolds} and \cref{thm:main_manifold}. This opens the way to new problems, wich are more naturally defined on standard manifolds like the $d$-dimensional sphere $S^d$ or the $d$-dimensional torus $\mathbb{T}^d \approx (S^1)^d$.

\subsection*{Background} The problem of decomposing $C^p$ functions as sums of squares has appeared in the context of symbolic calculus, in the proof of the Fefferman-Phong inequality, which is an important regularity result for partial differential operators (see \cite{Fefferman4673} for the original article, and \cite{bony1998} for the link with sum of squares decompositions, as well as \cite{tataru2002}). 
In this context, the following result is proved (with $C^{k,1}_{loc}$ denoting the set of $k$ times differentiable functions with locally Lispchitz $k$-th derivative): 

\begin{theorem}[Fefferman-Phong \cite{Fefferman4673}, Theorem 1.1 of \cite{bony2006nonnegative}]\label{thm:fefferman}
Let $\Omega$ be an open set of $\R^d$, $d \geq 1$ and $f \in C^{3,1}_{loc}(\Omega)$ be a non-negative function. Then $f$ can be written as a finite sum of squares of $C^{1,1}_{loc}(\Omega)$ functions. 
\end{theorem}

In the context of preserving regularity, a natural question which arises is whether increasing the regularity of $f$ can increase the regularity of the functions in a sum of square decomposition.
In \cite{bony2006nonnegative,bony2005}, it is shown that the general answer (under no further assumptions) is negative. More precisely, if $f$ is a function defined on a neighborhood of $0$, a local decomposition of $f$ around $0$ of class $\cc$ is a finite family $(f_i)_{ i \in I}$ of functions of class $\cc$ defined on an open neighborhood $U$ of $0$ such that $\sum_{i\in I}{f_i^2} = f$ on $U$.

\begin{theorem}[Theorem 2.1 of \cite{bony2006nonnegative}]\label{thm:impossibility}
In all the following cases, there exists $f \in C^{\infty}$ defined on an open neighborhood of $0$ in $\R^d$ such that the following holds: 
\begin{itemize}
    \item if $d \geq 4$, $f$ has no local decomposition of class $C^2$;
    \item if $d = 3$, $f$ has no local decomposition of class $C^3$.
\end{itemize}
\end{theorem}

The case $d = 1$ is explored in \cite{bony2005}: it is shown in Theorem 1 that if $f$ is of class $C^{2m}$ for $m$ finite, then $f$ can be written as the sum of squares of two functions of class $C^m$. Moreover, this is shown to be tight: there exists a function $f \in C^{2m}$ with no local decomposition as a sum of squares of functions of class $C^{m+k}$, for $k \geq 1$. The case $d = 2$ has been explored less in the literature (some results exist when dealing with flat minima, see for example Theorem 2 of \cite{bony2005}).

To summarize, these results show that without additional assumptions, as soon as the dimension is greater than $3$, inheriting the $C^p$ regularity properties of the function $f$ in the sum of squares decomposition is not possible in a satisfactory way, and motivates the introduction of additional geometric assumptions.

\paragraph{Polynomials} Decomposing non-negative polynomials as sums of squares has been related to important problems in algebraic geometry during the 20th century. In 1927, on his way to the resolution of Hilbert's 17th problem, Artin \cite{Artin1927} proved that any non-negative polynomial is a sum of squares of rational functions (that is formal fractions of polynomials $P(x)/Q(x)$). Moreover, Hilbert had earlier proved that there exist non-negative polynomials which cannot be written as sum of squares of polynomials in \cite{Hilbert1888} (for more than $3$ variables and with degree at least $6$ for example). In algebraic geometry, the set of SoS polynomials has very interesting properties, and finding sufficient conditions for a non-negative polynomial or even a positive polynomial to be a sum of squares is an important question. More generally, one usually wishes to understand under which sufficient conditions a polynomial $P$ which is non-negative (or positive) on an algebraic set, i.e., defined by polynomial inequalities of the form $Q_i \geq 0$ for polynomials $Q_i$, can be written in the form
$P=P_0 + \sum_{i=1}^N{P_i Q _i}$ where the $P_i$ are SoS.
The theoretical literature regroups these results under the name "Positivstellensatz". The most often seen in the SoS optmization literature are the Stengel \cite{Stengle1974}, Schm\"ugden \cite{schmugden1991} and Putinar \cite{Putinar1993} Positivstellens\"atzen.

If these algebraic geometry considerations seem far from applications and from decomposing smooth functions as sums of squares (indeed, polynomials are much more rigid than smooth functions) at first glance, they are actually related in two ways. 

    First, as smooth functions can be locally approximated by polynomials, results on polynomials give a good intuition of the difficulties one can encounter at the local level when decomposing a function as a sum of squares. Indeed, on the one hand, the general impossibility results proved in \cite{bony2005,bony2006nonnegative} (see \cref{thm:impossibility}) are obtained using Hilbert's theorem on the existence of non-negative polynomials which are not sum of squares. On the other hand, the fact that there is hope using our second-order assumptions is also due to the fact that second order non-negative polynomials can always be written as sums of squares.

Second, the certificates given by Positivstellensatz on the decomposability of certain non-negative polynomials can be algorithmically checked in some cases, using semi-definite programming. This has paved the way to so-called SoS hierarchies, and optimization of polynomial objective functions with polynomial constraints. These have been developed by Lasserre \cite{lasserre2010moments} (based on the Putinar Postivstellensatz \cite{Putinar1993}) and Parillo \cite{parillo2003} (based on the Stengel and Schm\"ugden Positivstellensatz \cite{Stengle1974,schmugden1991}). Using these theoretical results, they can provide certificates of lower bounds for certain optimization problems (or upper bound in the dual ``moment problem'', see \cite{lasserre2010moments}). Moreover, to have more interpretable results for these more applied settings, theses works have motivated more practical Positivstellensatz, like that in \cite{marshall2006}, which provides a condition for writing a polynomial with a finite set of zeros as as sum of squares (this condition is actually a second order condition which greatly resembles ours in the polynomial setting, although it deals more with the constraints $Q_i$).

\paragraph{$p$-times differentiable functions} In the same spirit as the polynomial hierarchies, recent works \cite{marteau2020non,rudi2020finding,rudi2021psd} have developed models and methods based on sum of squares of regular functions. The computational properties of these methods are based on the fact that regular functions can be well-approximated by functions of the form $\sum_{i}{\alpha_i k(\cdot,x_i)}$ where $k$ is a so-called positive definite kernel \cite{aronszajn1950theory} and can be adapted to the regularity.
In order to obtain guarantees on these methods, it is crucial to have the equivalent of Positivstellensatz in the case of regular functions. Contrary to the case of algebraic geometry, where such results existed for other purposes, there is a need to build such results for regular functions from scratch. Certain results like \cref{thm:simple_minima} have been presented. However, the aim of the present paper is to provide more general results, to be used in most situations.

\subsection*{Organisation of the work}

In \cref{sec:euclidean}, we formalize the different notions needed to state \cref{thm:general_Rd_intro} in the case of non-negative functions defined on open sets of $\R^d$. In particular, we start by presenting a local decomposition in \cref{thm:local_decomposition}, which will be the cornerstone of the work.
In \cref{sec:manifolds}, we extend \cref{thm:general_Rd_intro} to the manifold setting, and detail the procedure in which we glue local decompositions into a global one, using traditional tools from differential geometry.
In \cref{sec:proofs}, we formally prove \cref{thm:local_decomposition}. We finish by a discussion on the result presented in this paper, as well as possible extensions in \cref{sec:conclusion}.

\section{Decomposition as sums of squares given second order conditions (Euclidean case)}
\label{sec:euclidean}

In this section, we present our results on decomposing a $C^p$ function $f$ as a sum of squares of $C^{p-2}$ functions on open sets of $\R^d$.We start with a brief presentation of the notion of sub-manifold of $\R^d$ in \cref{secsec:submanifold}. It is the key geometric object we use to represent the set of zeros $\Z$ of the function $f$. In \cref{sec:local_sos_rd}, we present the cornerstone result of this paper in \cref{thm:local_decomposition}, as well as a sketch of its proof, which is done extensively in \cref{sec:proofs}. It shows that as soon as a non negative function has positive Hessian in the orthogonal direction to its zeros at a given point, then it can be decomposed as a sum of squares around that point. Finally, in \cref{secsec:global_sos_rd}, we present \cref{thm:main_Rd}, which shows that given a function defined on an open subset $\Omega$ of the Euclidean space $\R^d$, and under conditions on the Hessian of $f$ at its zeros, $f$ can be decomposed as a locally finite sum of squares of functions defined on $\Omega$.

\subsection*{Definitions and notations}

In general, given two topological sets $M$ and $N$ as well as $x_0 \in M$ and $y_0 \in N$, we will say that  $ \phi: (x_0,M) \rightarrow (y_0,N)$ is a local map satisfying a property $(P)$ if there exists an open neighborhood $U$ of $x_0$ in $M$ such that $\phi: U \rightarrow N$ is well defined, satisfies $\phi(x_0) = y_0$ and property $(P)$.
We will say that $\phi:U \subset \R^d \rightarrow \R^e$ defined on an open set $U$ is of class $C^k$ if it is $k$ times differentiable, and its derivatives of order $k$ are continuous.
For any function $\phi: (x,\R^d) \rightarrow \R^e$ of class $C^1$, we denote with $d \phi (x)$ its differential at $x$. It is an element of $\Hom(\R^d,\R^e)$ the set of linear maps from $\R^d$ to $\R^e$. We will write $d \phi(x) \xi$ or $d\phi(x)[\xi]$ the evaluation of $d\phi(x)$ at $\xi$. The Jacobian of $\phi$ at $x$ is the matrix $J_{\phi}(x) \in \R^{e \times d}$ which is the matrix of $d \phi(x)$ in the canonical bases. Writing the coordinates of $\phi$: $\phi = (\phi^1,...,\phi^e)$, we have $[J_{\phi}]_{ij} = \tfrac{\partial \phi^i}{\partial x^j}(x)$.

\subsection{Sub-manifolds of $\R^d$}\label{secsec:submanifold}
One of the main assumptions in order to achieve our results will be that the set of zeros of the non-negative function $f$ is a sub-manifold of $\R^d$. In this section, we restrict ourselves to introducing the definitions and results needed to state and prove those in this paper. For a more comprehensive introduction, see chapter 1 of \cite{lafontaine2015introduction}, section 2.2 of \cite{paulingeodiff} (in French) or \cite{spivak}. The notion of sub-manifold generalizes the notion of a curve in $\R^d$ (a one dimensional manifold) or a surface in $\R^d$ (a two dimensional manifold). Intuitively, a sub-manifold $N$ is a subset of $\R^d$ such that at each point $x \in N$, $N$ ``looks like'' $\R^{d_0}$ where $d_0$ is the dimension of the sub-manifold at $x$ (one for a line, two for a surface,etc.). Another way to put this is that $N$ can be locally parametrized by $\R^{d_0}$. To formalize this, we need the following definitions. We fix a subset $N \subset \R^d$.

A map $\phi: U \rightarrow \R^d$ defined on an open neighborhood $U$ of $0$ in $\R^{d_0}$ is said to be a local parameterization of $N$ around $x_0$ of class $C^k$ for $k\geq 1$ if $\phi$ is of class $C^k$, and if there exists an open set $V \subset \R^d$ such that the following conditions are satisfied: 
\begin{itemize}
    \item[(i)] $\phi(0) = x_0$, $\phi(U) = N\cap V$, and $\phi: U \rightarrow \phi(U)$ is a homeomorphism, i.e., it is bijective and has continuous inverse;
    \item[(ii)] its differential at $0$ is injective (one to one), i.e., $d \phi(t_0) \in \Hom(\R^{d_0},\R^d)$ is injective.
\end{itemize}

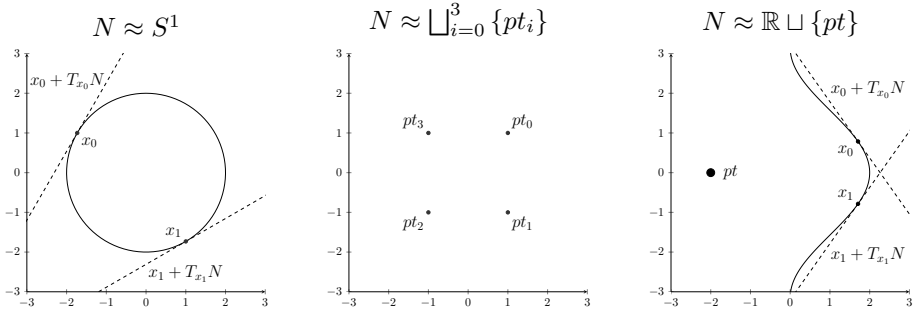
\begin{figure}\label{fig:sub_r2}
\centering 
\begin{tikzpicture}[scale = 0.45]  
\pgfmathsetmacro{\xz}{-pi/3};
\pgfmathsetmacro{\yz}{5*pi/6};
\pgfmathsetmacro{\r}{2};
\begin{axis}[
xmin=-3,xmax=3,
ymin=-3,ymax=3,
axis x line = bottom,
axis y line = left,
axis equal image,
]
\addplot[
thick,
domain=0:2*pi,
samples=100,]
({\r*sin(deg(x))},{\r*cos(deg(x))});
\addplot[
thick,
dashed,
domain=-3:3,
samples=100,]
({\r*sin(deg(\xz)) + x*\r*cos(deg(\xz))},{\r*cos(deg(\xz))-x*\r*sin(deg(\xz))});
    \node[black,below right] at 
    (axis cs:{\r*sin(deg(\xz))},{\r*cos(deg(\xz))})
    {\Large $x_0$}; 
    \node at 
    (axis cs:{\r*sin(deg(\xz))},{\r*cos(deg(\xz))})
    [circle, scale=0.3, draw=black!80,fill=black!80] {};
    \node[black,above] at 
    (axis cs:{\r*sin(deg(\xz))-0.25},{\r*cos(deg(\xz))+1})
    {\Large $x_0+T_{x_0}N$}; 
\addplot[
thick,
dashed,
domain=-3:3,
samples=100,]
({\r*sin(deg(\yz)) + x*\r*cos(deg(\yz))},{\r*cos(deg(\yz))-x*\r*sin(deg(\yz))});
    \node[black,above left] at 
    (axis cs:{\r*sin(deg(\yz))},{\r*cos(deg(\yz))})
    {\Large $x_1$}; 
    \node at 
    (axis cs:{\r*sin(deg(\yz))},{\r*cos(deg(\yz))})
    [circle, scale=0.3, draw=black!80,fill=black!80] {};
    \node[black,below] at 
    (axis cs:{\r*sin(deg(\yz))},{\r*cos(deg(\yz))-0.5})
    {\Large $x_1+T_{x_1}N$}; 
\end{axis}
\node[above] at (current bounding box.north) {$N \approx S^1$};
\end{tikzpicture}
\hspace{0.5cm}
\begin{tikzpicture}[scale = 0.45]  
\begin{axis}[
xmin=-3,xmax=3,
ymin=-3,ymax=3,
axis equal image, 
axis x line = bottom,
axis y line = left,
]
    \node[black,above right] at (axis cs:1,1) {\Large $pt_0$}; 
  \node at (axis cs:1,1) [circle, scale=0.3, draw=black!80,fill=black!80] {};
      \node[black,below right] at (axis cs:1,-1) {\Large $pt_1$}; 
  \node at (axis cs:1,-1) [circle, scale=0.3, draw=black!80,fill=black!80] {};
      \node[black,below left] at (axis cs:-1,-1) {\Large $pt_2$}; 
  \node at (axis cs:-1,-1) [circle, scale=0.3, draw=black!80,fill=black!80] {};
      \node[black,above left] at (axis cs:-1,1) {\Large $pt_3$}; 
  \node at (axis cs:-1,1) [circle, scale=0.3, draw=black!80,fill=black!80] {};
\end{axis}
\node[above] at (current bounding box.north) {$N \approx \bigsqcup_{i=0}^3{\{pt_i\}}$};
\end{tikzpicture}
\hspace{0.5cm}
\begin{tikzpicture}[scale = 0.45]  
\pgfmathsetmacro{\xzz}{-pi/4};
\pgfmathsetmacro{\yzz}{pi/4};
\begin{axis}[
xmin=-3,xmax=3,
ymin=-3,ymax=3,
axis equal image, 
axis x line = bottom,
axis y line = left,
]
\addplot[
    thick,
    domain=-pi:pi,
    samples = 100,
]
({1+cos(180*x/pi)},
{x});
\addplot[
    thick,
    dashed,
    domain=-pi:pi,
    samples = 100,
]
({1+cos(deg(\xzz)) -x*sin(deg(\xzz))},
{\xzz + x});
\addplot[
    thick,
    dashed,
    domain=-pi:pi,
    samples = 100,
]
({1+cos(deg(\yzz)) -x*sin(deg(\yzz))},
{\yzz + x});
\node at (axis cs:{1+cos(deg(\yzz))},{\yzz }) [circle, scale=0.3, draw=black!100,fill=black!100] {};
\node at (axis cs:{1+cos(deg(\xzz))},
{\xzz }) [circle, scale=0.3, draw=black!100,fill=black!100] {};
\node[black, below left] at (axis cs:{1+cos(deg(\yzz))},{\yzz })  {\Large $x_0$};
\node[black,above left] at (axis cs:{1+cos(deg(\xzz))},
{\xzz }) {\Large $x_1$};
\node[black, above] at (axis cs:{1.25+cos(deg(\yzz))},{1+\yzz })  {\Large $x_0+T_{x_0}N$};
\node[black,below] at (axis cs:{1.25+cos(deg(\xzz))},
{\xzz -1}) {\Large $x_1+T_{x_1}N$};
\node at (axis cs:-2,0) [circle, scale=0.7, draw=black!100,fill=black!100] {};
\node[black,right] at (axis cs:-2,0) {\Large $~~pt$};
\end{axis}
\node[above] at (current bounding box.north) {$N \approx \R \sqcup \{pt\}$};
\end{tikzpicture}

\caption{Examples of sub-manifolds of $\R^2$; points are denoted with $pt$. \textbf{Left:} connected sub-manifold of dimension $1$ (a circle). \textbf{Center:} a sub-manifold of $4$ connected components which are all points, i.e., of dimension $0$ (their tangent space is not represented since it is reduced to $\{0\}$). \textbf{Right:} a sub-manifold of two connected components, one point $pt$ of dimension $0$ and one of dimension $1$.}
\end{figure}

\begin{figure}\label{fig:sub_r3}
\centering
\begin{tikzpicture}[scale = 1.2,
declare function = {spherepar1(\x,\y) = 
{sin(deg(\y))*cos(deg(\x))}; 
},
declare function = {
spherepar2(\x,\y) = {sin(deg(\y))*sin(deg(\x))};
},
declare function = {
spherepar3(\x,\y) ={cos(deg(\y))};
},
declare function = {
der11(\x,\y) ={sin(deg(\x))};
},
declare function = {
der12(\x,\y) ={cos(deg(\x))};
},
declare function = {
der13(\x,\y) ={0};
},
declare function = {
der21(\x,\y) ={cos(deg(\x))*cos(deg(\y))*(1-2*sin(deg(\x))^2)/sqrt(1 + 4*cos(deg(\y))^2*cos(deg(\x))^4 -4*cos(deg(\y))^2*cos(deg(\x))^2)};
},
declare function = {
der22(\x,\y) ={sin(deg(\x))*cos(deg(\y))*(1-2*cos(deg(\x))^2)/sqrt(1 + 4*cos(deg(\y))^2*cos(deg(\x))^4 -4*cos(deg(\y))^2*cos(deg(\x))^2)};
},
declare function = {
der23(\x,\y) ={-sin(deg(\y))/sqrt(1 + 4*cos(deg(\y))^2*cos(deg(\x))^4 -4*cos(deg(\y))^2*cos(deg(\x))^2)};
},
]
\tikzstyle{every node}=[font=\small]
\begin{axis}[axis equal image,ymin = -1.25,ymax = 1.25,zmin = -1.25,zmax = 1.25,xmin=-1.25,xmax=1.25,ticklabel style = {font = \tiny}]
\pgfmathsetmacro{\a}{0}
\pgfmathsetmacro{\b}{pi/3}

\addplot3[surf,
      samples=20,
      domain=0:2*pi,y domain=0:pi,
      z buffer=sort,
      opacity = 0.4,
      fill = black!40!white,
      colorbar,
      colormap ={grayscale}{gray(0cm)=(0.5); gray(1cm)=(0.5)},
      ]
      ({spherepar1(x,y)},{spherepar2(x,y)},{spherepar3(x,y)});
\addplot3[surf,
      samples=10,
      domain=-1:1,y domain=-1:1,
      z buffer=sort,
      opacity = 0.4,
      fill = blue!40!white,
      colorbar,
      colormap ={grayscale}{gray(0cm)=(0.5); gray(1cm)=(0.5)},
      ]
      ({spherepar1(\a,\b) +x*der11(\a,\b) +y*der21(\a,\b)}, 
      {spherepar2(\a,\b) +x*der12(\a,\b) +y*der22(\a,\b)},
        {spherepar3(\a,\b) +x*der13(\a,\b) +y*der23(\a,\b)}
        );
        \node at (axis cs:{spherepar1(\a,\b)}, {spherepar2(\a,\b)},{spherepar3(\a,\b)}) [circle, scale=0.3, draw=black!100,fill=black!100] {};
        \node at (axis cs:{spherepar1(\a,\b)}, {spherepar2(\a,\b)},{spherepar3(\a,\b)}) [above] {$x_0$};
        \node at (axis cs:-0.3,0.75,1.1)  {$x_0+T_{x_0}N$};
\end{axis}
\end{tikzpicture}
\begin{tikzpicture}[scale = 1.2,
declare function = {tor1(\r,\s,\x,\y) = 
{-\s*(1+\r*cos(deg(\x)))*sin(deg(\y))}; 
},
declare function = {
tor2(\r,\s,\x,\y) = {\r*sin(deg(\x))};
},
declare function = {
tor3(\r,\s,\x,\y) ={\s*(1+\r*cos(deg(\x)))*cos(deg(\y))};
},
declare function = {
der11(\x,\y) ={cos(deg(\y))};
},
declare function = {
der12(\x,\y) ={0};
},
declare function = {
der13(\x,\y) ={sin(deg(\y))};
},
declare function = {
der21(\s,\x,\y) ={\s*sin(deg(\y))*sin(deg(\x))/sqrt(\s^2 + (1-\s^2)*cos(deg(\x))^2)};
},
declare function = {
der22(\s,\x,\y) ={cos(deg(\x))/sqrt(\s^2 + (1-\s^2)*cos(deg(\x))^2)};
},
declare function = {
der23(\s,\x,\y) ={-\s*cos(deg(\y))*sin(deg(\x))/sqrt(\s^2 + (1-\s^2)*cos(deg(\x))^2)};
},
]
\tikzstyle{every node}=[font=\small]
\pgfmathsetmacro{\r}{0.4}
\pgfmathsetmacro{\s}{0.7}
\pgfmathsetmacro{\a}{0}
\pgfmathsetmacro{\b}{0}
\pgfmathsetmacro{\la}{1}
\begin{axis}[axis equal image,ymin = -1.25,ymax = 1.25,zmin = -1.25,zmax = 1.25,xmin=-1.25,xmax=1.25,ticklabel style = {font = \tiny}]
\addplot3[
      samples=2,line width = 2pt,
      domain=0:5,opacity = 0.4,
      color = blue!60!white]
      ({0}, 
      {x},
        {0}
        );
\addplot3[
      variable = t,
      domain=0:1,samples y = 0,samples = 10]
      ({\la *t*t}, {t},{0} );
\addplot3[surf,
      samples=20,
      domain=0:2*pi,y domain=0:2*pi,
      z buffer=sort,
      opacity = 0.4,
      fill = black!40!white,
      colorbar,
      colormap ={grayscale}{gray(0cm)=(0.5); gray(1cm)=(0.5)},
      ]
      ({tor1(\r,\s,x,y)}, 
      {tor2(\r,\s,x,y)},
        {tor3(\r,\s,x,y)}
        );
        \addplot3[surf,
      samples=10,
      domain=-1:1,y domain=-1:1,
      z buffer=sort,
      opacity = 0.4,
      fill = blue!40!white,
      colorbar,
      colormap ={grayscale}{gray(0cm)=(0.5); gray(1cm)=(0.5)},
      ]
      ({tor1(\r,\s,\a,\b) +x*der11(\a,\b) +y*der21(\s,\a,\b)}, 
      {tor2(\r,\s,\a,\b) +x*der12(\a,\b) +y*der22(\s,\a,\b)},
        {tor3(\r,\s,\a,\b) +x*der13(\a,\b) +y*der23(\s,\a,\b)}
        );
        \node at (axis cs:{tor1(\r,\s,\a,\b)}, {tor2(\r,\s,\a,\b)},{tor3(\r,\s,\a,\b)}) [circle, scale=0.3, draw=black!100,fill=black!100] {};
        \node at (axis cs:{tor1(\r,\s,\a,\b)}, {tor2(\r,\s,\a,\b)},{tor3(\r,\s,\a,\b)}) [above] {$x_0$};
        \addplot3[
      samples=2,
      domain=-2:0,line width =2pt,
      color = blue!40!white,opacity = 0.5]
      ({0}, 
      {x},
        {0}
        );
        \addplot3[
      samples=10, samples y=0,
      domain=-1:0,thick,
      ]
      ({\la*x^2}, 
      {x},
        {0}
        );
        \node at (axis cs:0,0,0) [circle, scale=0.3, draw=black!100,fill=black!100] {};
        \node at (axis cs:-0.15,0,0.15) {$x_1$};
        \node at (axis cs:0,-1,-0.8) {$x_1 + T_{x_1}N$};
        \node at (axis cs:-0.4,1,1) {$x_0 + T_{x_0}N$};
\end{axis}
\end{tikzpicture}

\vspace*{-.2cm}

\caption{Two examples of sub-manifolds of $\R^3$. The blue affine spaces represent tangent spaces. \textbf{Left:} connected sub-manifold of dimension $2$ (the sphere $S^1$). \textbf{Right:} a sub-manifold of two connected components, one of dimension $2$ (homeomorphic to the torus $\mathbb{T}^2$ on which lies $x_0$),  and one of dimension $1$ on which lies $x_1$.}
\end{figure}
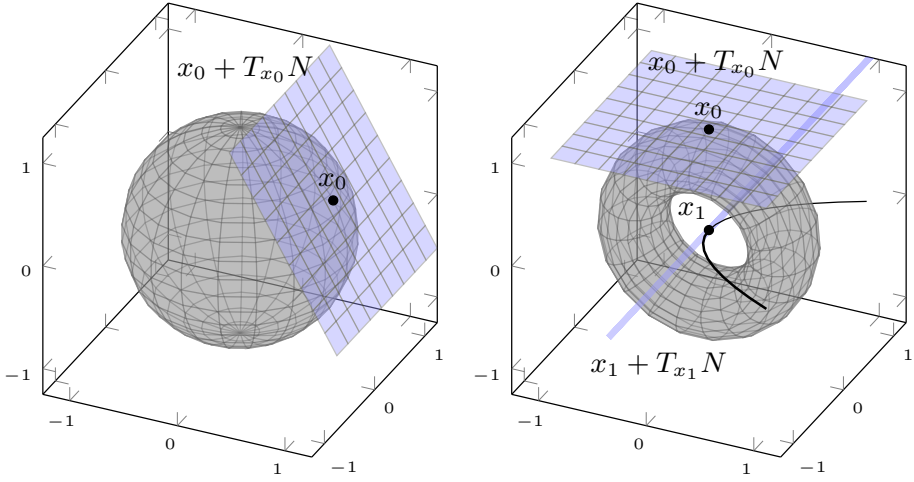

The second condition guarantees that the local dimension of $N$ is indeed $d_0$, that $\phi$ is not an over-parameterization. $N$ is said to be a sub-manifold of $\R^d$ and of class $C^k$ if there exists a local parameterization $\phi$ of class $C^k$ around each point $x \in N$. Given a point $x \in N$, the dimension $d_x$ of the local parametrization is independent of the parametrization (two local parametrizations will necessarily be of same dimension); it is called the dimension of $N$ at $x$.  Similarly, the subspace $T_{x}N:= d\phi(x)\R^{d_x}$, which is a subspace of $\R^d$ of dimension $d_x$ is independent of the local parametrization: it is the linear approximation of $N$ at $x$ and is called the tangent space to $N$ at $x$ (see \cref{fig:sub_r2} and \cref{fig:sub_r3} for more visual representations). 

A sub-manifold $N$ of $\R^d$ is said to be connected if it cannot be written as a union of disjoint open sets. Equivalently, it is connected if any two points in $N$ can be connected by a continuous path $\gamma: [0,1]\rightarrow N$. On a connected sub-manifold $N$, the dimension $d_x$ is the same at every point $x$, it is called the dimension of the connected sub-manifold $N$. This implies that all the tangent spaces $T_x N$ have the same dimension.

\begin{example}
All open sets of $\R^d$ are sub-manifolds of $\R^d$. The $d$-dimensional sphere $S^d$ is a sub-manifold of $\R^{d+1}$. $S^1$ is represented in the left hand side (l.h.s.) of \cref{fig:sub_r2} and $S^2$ in the l.h.s. of \cref{fig:sub_r3}. Given a sub-manifold $N$ of $R^d$, the intersection of $N$ with any open set of $\R^d$ is a sub-manifold of $\R^d$.
\end{example}

If $U_i$ is a family of disjoint open sets each containing a connected sub-manifold $N_i$ of $\R^d$, it is clear the the disjoint union $\sqcup_{i \in I}{N_i}$ is also a sub-manifold. Conversely, any sub-manifold can be decomposed into its connected components $N_i$; moreover, one can find a family of disjoint open sets $U_i$ such that $N_i \subset U_i$ (see \cref{lm:open_sets}).

These results, their proof and their broader context can be found in \cite{lafontaine2015introduction} in chapter 1.5. In particular, Theorem 1.21 presents equivalent definitions of a sub-manifold. Section 2.2 of \cite{paulingeodiff} is also a good reference (in French).

\subsection{Local decomposition as a sum of squares}\label{sec:local_sos_rd}

In this section, $f$ will always denote a non-negative function defined on an open set of $\R^d$. We will assume that $f$ is of class $C^p$ for $p \geq 2$. We will also denote with $\Z$ the set of zeros of $f$, i.e., the set of zeros of $f$. In this section, we will make local assumptions on the Hessian of $f$ at points $x \in \Z$ such that the function $f$ can be decomposed as a sum of squares locally around $x$.

We will denote with $d^2 f(x)$ the second differential of $f$ (which we will sometimes call abusively its Hessian), which is a symmetric bilinear form on $\R^d$. We denote with $d^2f(x)[\xi,\eta]$ its evaluation on vectors $\xi,\eta$. We denote with $\nabla^2 f(x) \in \R^{d\times d}$ the Hessian matrix of $f$ at $x$, which is the matrix of $d^2f(x)$ in the canonical basis of $\R^d$, and we have $d^2 f(x)[\xi,\eta] = \eta^{\top}\nabla^2 f(x)\xi$. For any vector sub-space space $S \subset \R^d$, and any bilinear form $H$ on $\R^d$, we denote with $H|_{S}$ the restriction of $H$ to $S$, which is a bilinear form on $F$. We say that a bilinear form $H$ is positive semi-definite if $H[\xi,\xi] \geq 0$ for any $\xi \in \R^d$, and is positive definite if $H[\xi,\xi] >0$ for all $\xi \in \R^d \setminus{\{0\}}$. We use the same terminology for matrices.

    
We are now ready to state \cref{thm:local_decomposition}, which is the cornerstone of this work. For the rest of this section (\cref{sec:local_sos_rd}), let $x_0 \in \R^d$ and $f: (x_0,\R^d) \rightarrow \R$ be a non-negative $C^p$ function, for $p\geq 2$, such that $f(x_0) = 0$. We claim that if there is a sub-manifold of class $C^1$ and of dimension $d_0$ around $x_0$ of zeros of $f$, and if the Hessian of $f$ at $x_0$ has rank $d-d_0$ (which we will call the normal Hessian condition), then it can be decomposed as a sum of squares as in \cref{eq:sos_form}. 

\begin{definition}[normal Hessian condition]\label{df:nhc_rd}
 Let $\Z$ denote the set of zeros of $f$. We say that $f$ satisfies the {\em normal Hessian condition} (NHC) at $x_0$ if there exists a dimension $0 \leq d_0 \leq d$ and a sub-manifold $N$ of class $C^k$ with $k \geq 1$ and of dimension $d_0$ such that $x_0 \in N \subset\Z$, and on one of the following equivalent conditions is satisfied: 

\begin{itemize}
    \item[(i)] the rank of $\nabla^2 f(x_0)$ at $x_0$ is $d-d_0$;
    \item[(ii)] the restriction of $d^2 f(x_0)$ to $T_{x_0}N^{\perp}$ is positive definite.
\end{itemize}
\end{definition}

\begin{figure}\label{fig:NHC_proof}
\centering
\begin{tikzpicture}[scale = 0.65]
\tikzstyle{every node}=[font=\Large]
\begin{axis}[xlabel = $x_1$,ylabel=$x_2$]
\addplot3[
    surf,
    opacity = 0.4,
    domain=-4:4,
    y domain=-4:4
]
{(x - sin(180*y/pi)/2)^2};
\addplot3[
    line width = 3pt,
    domain=-4:4,
    samples = 100,
    samples y=0,
]
({sin(180*x/pi)/2},
{x},
{0});
\addplot3[
    line width = 3pt,
    domain=-4:4,
    samples = 100,
    samples y=0,
]
({sin(180*x/pi)/2},
{x},
{0});
\node[above left] at (axis cs:{-0.4,0,0}) { $x_0$};
\node at (axis cs:0,0,0) [circle, scale=0.7, draw=black!100,fill=black!100] {};
\end{axis}
\end{tikzpicture}
\begin{tikzpicture}[scale = 0.65]
\tikzstyle{every node}=[font=\Large]
\begin{axis}[axis equal image,
xmin=-4,xmax=4,
ymin=-4,ymax=4,
xlabel = $x_1$,ylabel=$x_2$,zlabel=$f(x)$,
]
\addplot[
thick,
domain=-4:4,
samples=100,]
({sin(180*x/pi)/2},{x});
\addplot[
dashed,
domain=-4:4,
samples=100,]
({x/2},{x});
\addplot[
dotted,
domain=-4:4,
samples=100,]
({-x},{x/2});
\node[] at (axis cs:{2,2}) {$T_{x_0}N$};
\node[right] at (axis cs:{0,-2.5}) {$N$};
\node[above right] at (axis cs:{-2,1}) {$T_{x_0}N^{\perp}$};
\node[below left] at (axis cs:{-0.4,0}) {$x_0$};
\node at (axis cs:0,0) [circle, scale=0.7, draw=black!100,fill=black!100] {};
  \draw[very thick,->,red] (axis cs:{0,0}) -- (axis cs:{-0.89,0.44});
\end{axis}
\end{tikzpicture}
\begin{tikzpicture}[scale = 0.65]
\tikzstyle{every node}=[font=\Large]
\begin{axis}[axis equal image,
xmin=-4,xmax=4,
ymin=-4,ymax=4,
xlabel=$\xpar$,
ylabel=$\xperp$,
]
\def\r{2.}
\def\s{1.}
\addplot[
thick,
domain=-6:6,
samples=100,]
({-0.44*sin(180*x/pi)/2-0.89*x},{-0.44*x +  0.89*sin(180*x/pi)/2});
\addplot[
dashed,
domain=-4:4,
samples=100,]
({x},{0});
\addplot[
dotted,
domain=-4:4,
samples=100,]
({0},{x});
\node[above] at (axis cs:{-3,0}) {$T_{x_0}N$};
\node[right] at (axis cs:{-3,-2.5}) {$N$};
\node[right] at (axis cs:{0,3}) {$T_{x_0}N^{\perp}$};
\node[below left] at (axis cs:{-0,0}) {$x_0$};
\node at (axis cs:0,0) [circle, scale=0.7, draw=black!100,fill=black!100] {};
\draw[->] (axis cs:2.5,0) -- node[right]{$\varphi(\xpar)$} (axis cs:2.5,0.8);
\draw[very thick,->,red](axis cs:0,0)--(axis cs:0,1);
\end{axis}
\end{tikzpicture}
\begin{tikzpicture}[scale = 0.95]
        \begin{axis}[
            xmin = -4, xmax = 4,
            ymin = -4, ymax = 4,
            axis equal image,
            view = {0}{90},
            scale = 0.6,
            title = {Vector field $\xi(\xperp,\xpar) \in T_{x_0}N^{\perp}$},
            xlabel = {$\xpar$},
            ylabel = {$\xperp$},
            colorbar,
            colorbar style = {
                ylabel = {Vector Length}
            }
        ]
            \addplot3[
                point meta = {sqrt((-0.44*x+0.89*y - sin(-0.44*y-0.89*x)/2)^2)},
                quiver = {
                    u = {0},
                    v = {-0.44*x+0.89*y - sin(-0.44*y-0.89*x)/2},
                    scale arrows = 0.25,
                },
                quiver/colored = {mapped color},
                -stealth,
                domain = -4:4,
                domain y = -4:4,
            ] {0};
        \end{axis}
\end{tikzpicture}

\vspace*{-.2cm}

\caption{Local view of the function around a minimum lying on a $1$-dimensional manifold. \textbf{Top left:} function around the minimum $x_0$. \textbf{Top right:} decomposition of $\R^2$ at $x_0$ between tangent space and normal tangent space $T_{x_0} N + T_{x_0}N^{\perp}$, and positive eigen-vector of the Hessian in red. \textbf{Bottom left:} reparametrization in the right coordinate system, and representation of the map $\varphi(x)$ given by the Morse Lemma. \textbf{Bottom right:} vector field $\xi$ given by the Morse lemma.}
\end{figure}
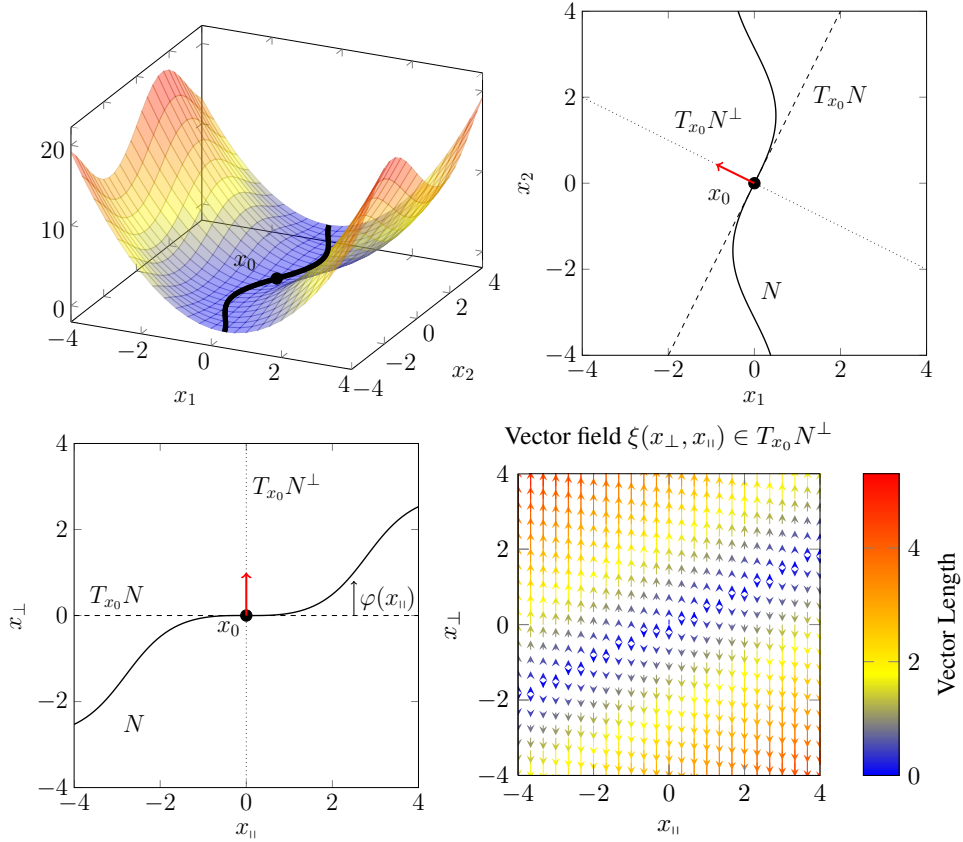

The complete proof of the equivalence of these conditions as well as the proof of \cref{thm:local_decomposition} can be found in \cref{sec:proof_local_decomposition}. To illustrate the definition of the normal Hessian condition, we refer to 
\cref{fig:NHC_proof}
which represents the local behavior of functions $f$ defined locally around a point $x_0 \in \R^2$ in the set of zeros and which satisfies the NHC for $d_0=1$.

\begin{theorem}\label{thm:local_decomposition}
If $f$ satisfies the NHC at $x_0$ (\cref{df:nhc_rd}) with regularity $k$ and dimension $d_0$, there exists an open neighborhood $U$ of $x_0$ in $\R^d$ on which $f$ is defined and such that $U \cap \Z$ is a sub-manifold of $\R^d$ of dimension $d_0$ and of class $C^{\max(k,p-1)}$, and there exist functions $f_i \in C^{p-2}(U)$ where $1 \leq i \leq d-d_0$ such that 
\begin{equation}
    \label{eq:local_form}
    \forall x \in U,~ f(x) = \sum_{i=1}^{d-d_0}{f_i^2(x)}.
\end{equation}
\end{theorem}

\begin{proof}[Main steps of the proof]
The main steps of this proof are represented geometrically in \cref{fig:NHC_proof}.

\textit{Step 1.} We show that under the NHC at $x_0$, we have $T_{x_0}N = \ker(\nabla^2f(x_0))$ and hence that $d^2 f(x_0)|_{T_{x_0}N^\perp}$ is positive definite. 

\textit{Step 2.} Re-parametrizing $f$ on a basis adapted to $T_{x_0}N^{\perp} \oplus T_{x_0}N$ as $f(\xperp,\xpar)$, we apply the Morse lemma (see \cref{lm:morse_hormander}), which decomposes the function $f$ in the form 
\begin{equation} f(\xperp,\xpar) = f(\varphi(\xpar),\xpar) + \tfrac{1}{2}d^2 
\label{eq:res_morse_lemma}
f(x_0)|_{T_{x_0}N^\perp}[\xi(\xperp,\xpar),\xi(\xperp,\xpar)],
\end{equation}
for a certain function $\varphi$ of class $C^{p-1}$ and $\xi$ of class $C^{p-2}$ in a certain open set around $x_0$ (for an easy visualization, see \cref{fig:NHC_proof}).

\textit{Step 3.} We show that the second term of the right hand side of \cref{eq:res_morse_lemma} can actually be seen as a sum of squares of $d-d_0$ functions of class $C^{p-2}$.

\textit{Step 4.} We characterize the manifold of zeros around $x_0$.

\textit{Step 5.} We show that the first term of the result of the Morse lemma is equal to zero using the previous characterization, which shows \cref{eq:local_form}.
\end{proof}

\begin{example}[case where $d_0 = 0$]
When $d_0 = 0$, the NHC at $x_0$ is simply the SHC \cref{eq:sc_1}, that is the condition that $x_0$ be a strict minimum. In that case, \cref{thm:local_decomposition} simply states that there exists an open neighborhood $U$ of $x_0$ such that $U\cap \Z = \{x_0\}$ and on which $f$ can be decomposed as the sum of $d$ squares.
\end{example}

\begin{remark}[Smoothing effect]
Note that \cref{thm:local_decomposition} induces a smoothing effect: indeed, if we simply assume that there exists a $d_0$ dimensional manifold of class $C^1$ of zeros satisfying the NHC, one sees that this manifold is actually of class $C^{p-1}$ in a neighborhood of~$x_0$. 
\end{remark}

\subsection{Global decomposition as a sum of squares for functions on $\R^d$}\label{secsec:global_sos_rd}

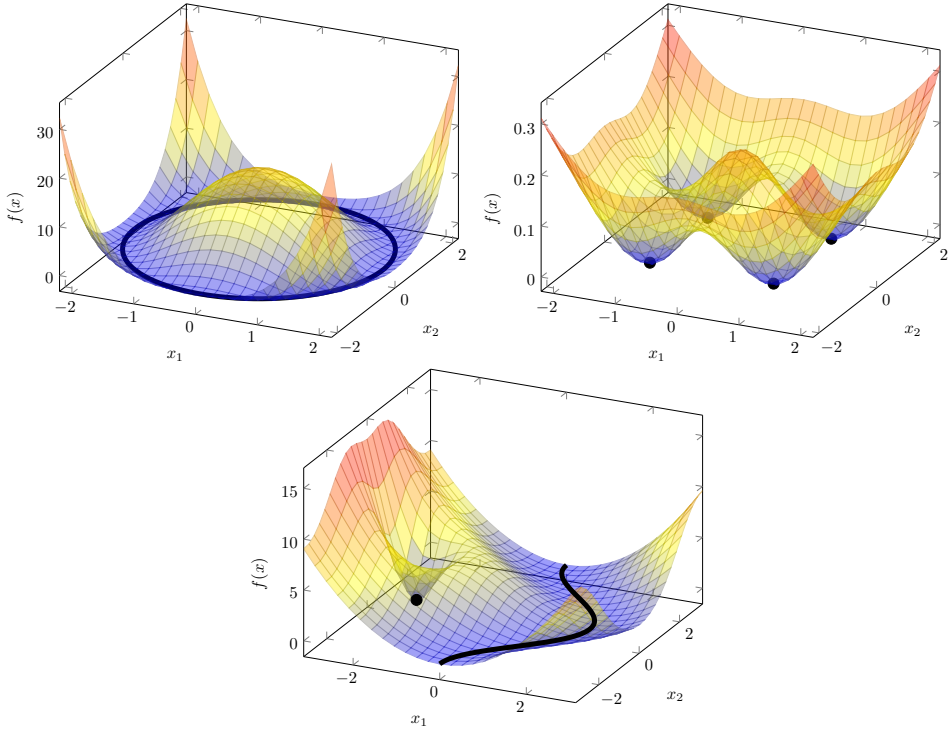
\begin{figure}\label{fig:examples_ghc}
\centering
\begin{tikzpicture}[scale = 0.63]
\begin{axis}[xlabel = $x_1$,ylabel=$x_2$,zlabel=$f(x)$]
\addplot3[
    line width = 3pt,
    domain=-pi:pi,
    samples = 100,
    samples y=0,
    xlabel = $x$,
    ylabel = $y$,
    zlabel = $z$,
]
({2*cos(180*x/pi)},
{2*sin(180*x/pi)},
{0});
\addplot3[
    opacity=0.4,
    surf,
    domain=-2.2:2.2,
    y domain=-2.2:2.2,
]
{(x^2 + y^2 -4)^2};
\end{axis}
\end{tikzpicture}
\begin{tikzpicture}[scale = 0.63]
\begin{axis}[xlabel = $x_1$,ylabel=$x_2$,zlabel=$f(x)$]
  \node at (axis cs:1,1,0) [circle, scale=0.7, draw=black!100,fill=black!100] {};
  \node at (axis cs:1,-1,0) [circle, scale=0.7, draw=black!100,fill=black!100] {};
  \node at (axis cs:-1,1,0) [circle, scale=0.7, draw=black!100,fill=black!100] {};
  \node at (axis cs:-1,-1,0) [circle, scale=0.7, draw=black!100,fill=black!100] {};
\addplot3[
    opacity=0.4,
    surf,
    domain=-2.2:2.2,
    y domain=-2.2:2.2,
]
{((x-1)^2 + (y-1)^2)*((x-1)^2 + (y+1)^2)*((x+1)^2 + (y-1)^2)*((x+1)^2 + (y+1)^2)*exp(-4 *ln(3+x^2 + y^2))};
\end{axis}
\end{tikzpicture}

\begin{tikzpicture}[scale = 0.63]
\begin{axis}[xlabel = $x_1$,ylabel=$x_2$,zlabel=$f(x)$]
\addplot3[
    surf,
    opacity = 0.4,
    domain=-pi:pi,
    y domain=-pi:pi
]
{((x - 1- cos(180*y/pi))^2)*((x+2)^2 + y^2)*exp(-ln(1+(x+2)^2 + y^2))};
\addplot3[
    line width = 3pt,
    domain=-pi:pi,
    samples = 100,
    samples y=0,
]
({1+cos(180*x/pi)},
{x},
{0});
\node at (axis cs:-2,0,0) [circle, scale=0.7, draw=black!100,fill=black!100] {};
\end{axis}
\end{tikzpicture}
\caption{Example of functions $f$ which satisfy the global normal Hessian condition, with sub-manifolds $\Z$ of zeros corresponding to the sub-manifolds presented in \cref{fig:sub_r2} in the same order}
\end{figure}

In this section, we fix $f$ to be a non-negative $C^p$ function defined on an open subset $\Omega$ of $\R^d$. Once again, we assume $p\geq 2$.
The goal is to find conditions on $f$ to be written as a sum of squares of functions defined on $\Omega$. These conditions will be that the NHC holds at every $x_0 \in \Z$.
We will start by reformulating this assumption in a more global and geometric way. We introduce the following definition of a manifold to which $f$ is positively normal.

\begin{definition}[positive normally to a sub-manifold]\label{df:positive_normal}
Let $N$ be a sub-manifold of $\R^d$ of class $C^k$ for $k\geq 1$ and included in $\Omega$. We say that $f$ is positive normally to $N$ if:
\begin{itemize}
    \item[a)] $N$ is included in the set of critical points of $f$ ($df(x) = 0$ for all $x \in N$);
    \item[b)] for any $x_0 \in N \cap \Omega$, if $d_0$ is the local dimension of $N$ at $x_0$, there exists a subspace $S \subset \R^d$ of dimension $d-d_0$ such that $d^2 f(x_0)|_{S}$ positive definite. 
\end{itemize}
\end{definition}

The intuition of this definition is that if $f$ is positively normal to $N$, then $f$ grows quadratically normally to $N$, which is a local minimum valley on which $f$ is constant. Note that there can be more than one connected component in $N$: this will correspond to multiple local minima valleys (see second and third examples in \cref{fig:examples_ghc}).
We now reformulate the fact that the NHC holds at every point in $\Z$ as a more geometric global assumption, using \cref{df:positive_normal}.

\begin{lemma}[Global normal Hessian condition]\label{df:global_nhc}
The following statements are equivalent and define the global NHC condition: 
\begin{enumerate}[label=\roman*]
\item[(i)] for all $x_0 \in \Z$, $f$ satisfies the NHC; 
\item[(ii)] $\Z$ is a sub-manifold of $\R^d$ (not necessarily connected) of class $C^1$ such that  the Hessian of $f$ is positive normally to $\Z$; 
\item[(iii)] $\Z$ is a sub-manifold of $\R^d$ (not necessarily connected) of class $C^{p-1}$ such that the Hessian of $f$ is positive normally to $\Z$. 
\end{enumerate}
\end{lemma}

This equivalence is a direct consequence of the local description of $\Z$ obtained in \cref{thm:local_decomposition} under the local NHC. 
Examples of functions satisfying the global normal Hessian condition can be found in \cref{fig:examples_ghc}. They have manifolds of zeros which are depicted in the same order in \cref{fig:sub_r2}.
Under this geometric condition, we will show in \cref{thm:main_Rd} that $f$ can be written as a sum of squares of $C^{p-2}$ functions with locally finite support, defined below. 

\noindent\textbf{Locally finite support.} Let $X$ be a topological space (see \cite{topologie} for full definitions). We say that a family $(S_i)$ of subsets of $X$ is locally finite if for every $x \in X$, there exists an open set $U_x$ containing $x$ which intersects a finite number of the $S_i$, i.e.,  $|\{i \in I ~:~ U_x \cap S_i \neq \emptyset\}| < \infty$. A family $(f_i)$ of functions on a topological space $X$ has locally finite support if the family of supports $(\supp(f_i))_{i \in I}$ is locally finite (recall that $\supp(f_i) = \overline{\{x~:~f(x)\neq 0\}})$. In particular, if $(f_i)$ has locally finite support, the function $\sum_{i \in I}{f_i^2}$ is well defined and it is also of class $C^{q}$ if the functions are of class $C^q$. Using this terminology, the global result can be stated as follows.

\begin{theorem}
\label{thm:main_Rd} 
If $f$ satisfies the global normal Hessian condition in  \cref{df:global_nhc}, there exists an at most countable family $(f_i)_{i \in I} \in (C^{p-2}(\Omega))^I$ with locally finite support such that  
\begin{equation}
    \label{eq:res_Rd}
    \forall x \in \Omega,~ f(x) = \sum_{i\in I}{f_i(x)^2}.
\end{equation}
Moreover: 
\begin{itemize}
    \item if $f$ satisfies the strict Hessian condition, $\Z$ is discrete and we can find such a decomposition such that $|I| \leq d+1$.
    \item if $\Z$ is compact, then $|I|$ can be taken to be finite.
    \end{itemize}
\end{theorem}

For the formal proof of this result, we refer to the next section, where this result will be proved more generally for functions defined on manifolds (see \cref{thm:main_manifold} and \cref{sec:manifolds}).

\begin{proof}[Main steps of the proof]
 For subtleties pertaining to the SHC case, we refer to \cref{sec:manifolds}. The gluing done in that section is slightly more elaborate.

\textit{Step 1.} Since the local NHC holds at any point in $\Z$, using \cref{thm:local_decomposition} shows that at any point~$x$, there exists an open neighborhood $U_x$ of $x$ , an integer $n_x$ and  functions $(f_{x,j})_{1 \leq j \leq n_x}$ of class $C^{p-2}$ on $U_x$ such that $f = \sum_{j=1}^{n_x}{f_{x,j}^2}$ on $U_x$. The collection of sets $U_x$ is then an open covering of $\Z$. Since $\R^d$ is Hausdorff and second-countable (see \cref{sec:manifolds} for precise definitions), only an at most countable subsets of them are necessary to  cover $\Z$ (even a finite number if $\Z$ is included in a compact, since it is then itself a compact). Denote with $(U_{i})_{i \in I}$ this open covering, and replace $x$ by $i$ to denote the associated $f_{i,j}$ and $n_i$. 

\textit{Step 2.} Since $\Z$ is closed, as the set of zeros of a continuous function, the set $U_{>0}:= \{x \in \Omega ~:~ f(x)>0\}$ is open and the map $f_1:= \sqrt{f}: U_{>0} \rightarrow \R$ is of class $C^p$ and satisfies $f_1^2 = f$. 
We can therefore add $U_{>0}$ to the collection $(U_i)$ and still guarantee the following property: for all $i \in I$, there exists $n_i \in \N$ and $f_{i,j} \in C^{p-2}(U_i)$ such that $f = \sum_{j=1}^{n_i}{f_{i,j}^2}$. Moreover, $(U_i)$ becomes an open covering of $\Omega$; in particular, if $U_i$ was a finite covering of $\Z$, it now becomes a finite covering of $\Omega$.

\textit{Step 3.} Using \cref{lm:glue_square}, we can take a partition of unity $(\chi_i)$ adapted to the open covering $\bigcup_{i} U_i$ such that $\sum_{i}{\chi_i^2} = 1$ and which is locally finite. Define $\tilde{f}_{i,j} = f_{i,j} \chi_{i}$ which is now defined on the whole of $\Omega$ (indeed, it can be extended as zero to $\Omega \setminus{U_i}$ since the support ot $\chi_i$ is included in $U_i$). The $\tilde{f}_{i,j}$ satisfy $\sum_{i,j}{\tilde{f}_{i,j}^2} = f$ on the whole of $\Omega$, and is a finite family if the covering $U_i$ is finite (if $\Z$ is assumed to be compact for example).
\end{proof}

\section{Global decomposition as a sum of squares for functions on manifolds}\label{sec:manifolds}

In this section, we present results analogous to those of \cref{sec:euclidean} but in the more general context of manifolds. After a brief recap on the terminology and definitions related to manifolds, in \cref{sec:reformulation_assumptions}, we will adapt the definitions of the local and global normal Hessian conditions, as well as state the equivalent result to \cref{thm:local_decomposition} in the context manifolds. In \cref{sec:gluing}, we will introduce the tools to glue local decompositions as sum of squares together. Finally, in \cref{sec:general_manifold}, we prove \cref{thm:main_manifold}, the equivalent of \cref{thm:main_Rd} in the broader context of manifolds.

\subsection*{Additional definitions and notations for manifolds}

In this section, we introduce the basic definitions we will need concerning manifold. For more formal introductions to manifolds, we refer to \cite{lafontaine2015introduction,paulingeodiff,spivak}. Informally, a manifold of dimension $d$ is a set which ``looks like $\R^d$'' locally. 
This means that at every point $x \in M$, we can find a chart $\phi$ which topologically maps a neighborhood $U$ of $x$ to an open set of $\R^d$.

More generally, we define a chart on a topological space $M$ as a map $\phi: U \rightarrow \R^d$ for some $d \in \N$, defined on an open set $U$ of $M$, and which is a homeomorphism onto its image.
We define a manifold $M$ as a second-countable\footnote{A topological space is said to be second countable if there exists a countable sequence of open sets $U_n$ such that any open set $U$ in the topology is a reunion of a part of the $U_n$.}, Hausdorff\footnote{A topological space is Hausdorff if for any two points $x \neq x^{\prime}$, there exists two open sets $U,V$ such that $x \in U$ and $x^{\prime} \in V$ and $U\cap V = \emptyset$} topological space equipped with a collection ${\cal A} = (\phi_i)_{i \in I}$ of charts such that 
\begin{itemize} 
\item[(i)]all transition maps $\phi_i \circ \phi_j^{-1}: \phi_j(U_j \cap U_i) \rightarrow \phi_i(U_i \cap U_j)$ are homeomorphisms;
\item[(ii)] the charts cover $M$ entirely, i.e. $M = \bigcup_{i \in I}{U_i}$.
\end{itemize}
The set ${\cal A}$ is called an atlas. The manifold $M$ is said to be of class $C^k$ for $k \geq 0$ if all the transition maps $\phi_i \circ \phi_j^{-1}$ are of class $C^k$. It is said to be of dimension $d$ if all its charts are in $\R^d$. As for sub-manifolds of $\R^n$, a manifold can always be decomposed as the union of its connected components, and the dimension is the same on each connected component. Note that with this definition, the restriction of a manifold to any open set is still a manifold (just by restricting the charts).

If $M$ is a manifold of class at least $C^1$, we can define at each point its tangent space $T_x M$. Informally, this set $T_x M$ is all the possible derivatives $\gamma^{\prime}(0)$ of curves $\gamma: I \rightarrow M$ defined on an open interval $I$ around $0$ such that $\gamma(0) = x$. Of course $\gamma^{\prime}(0)$ is not yet formally defined. 
Formally, $T_{x}M$ can be defined as the classes of $C^1$ curves defined on an open interval $I$ around $0$ such that $\gamma(0) = x$, where we identify two curves $\gamma$ and $\tilde{\gamma}$ if $(\phi \circ \gamma)^{\prime}(0) =(\phi \circ \tilde{\gamma})^{\prime}(0)$ for a (or equivalently any) chart $\phi$ of $M$ around $x$. We denote with $[\gamma] \in T_{x}M$ the equivalence class of a curve $\gamma$. It can be shown that $T_x M$ is a vector space (with the natural definition $\lambda [\gamma] + \mu[\tilde{\gamma}] = [\lambda \gamma + \mu \tilde{\gamma}]$) and that it is of dimension $d$ where $d$ is the dimension of $M$ at~$x$.

A map $g: M \rightarrow \R^p$ defined on a manifold $M$ is said to be of class $C^q$ if $M$ is of class at least $C^q$ and if for any chart $\phi$, the map $g \circ \phi^{-1}$ is of class $C^q$. If $q \geq 1$, then we can define the differential of $g$ at any point $x \in M$ as $df(x) [\gamma] = (f \circ \gamma)^{\prime}(0)$. Hence, $df(x) \in \Hom(T_x M,\R^p)$.

\begin{example}
All sub-manifolds of $\R^d$ are manifolds. The notions of regularity, dimension, and tangent space coincide.
\end{example}

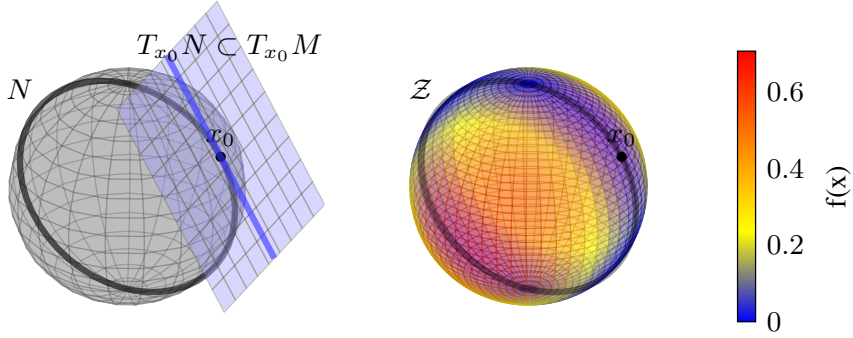
\begin{figure}\label{fig:manifold_fig}
\centering
\begin{tikzpicture}[scale = 1.2,
declare function = {spherepar1(\x,\y) = 
{sin(deg(\y))*cos(deg(\x))}; 
},
declare function = {
spherepar2(\x,\y) = {sin(deg(\y))*sin(deg(\x))};
},
declare function = {
spherepar3(\x,\y) ={cos(deg(\y))};
},
declare function = {
der11(\x,\y) ={sin(deg(\x))};
},
declare function = {
der12(\x,\y) ={cos(deg(\x))};
},
declare function = {
der13(\x,\y) ={0};
},
declare function = {
der21(\x,\y) ={cos(deg(\x))*cos(deg(\y))*(1-2*sin(deg(\x))^2)/sqrt(1 + 4*cos(deg(\y))^2*cos(deg(\x))^4 -4*cos(deg(\y))^2*cos(deg(\x))^2)};
},
declare function = {
der22(\x,\y) ={sin(deg(\x))*cos(deg(\y))*(1-2*cos(deg(\x))^2)/sqrt(1 + 4*cos(deg(\y))^2*cos(deg(\x))^4 -4*cos(deg(\y))^2*cos(deg(\x))^2)};
},
declare function = {
der23(\x,\y) ={-sin(deg(\y))/sqrt(1 + 4*cos(deg(\y))^2*cos(deg(\x))^4 -4*cos(deg(\y))^2*cos(deg(\x))^2)};
},
]
\tikzstyle{every node}=[font=\small]
\begin{axis}[axis equal image,ymin = -1.25,ymax = 1.25,zmin = -1.25,zmax = 1.25,xmin=-1.25,xmax=1.25,ticks=none,axis line style={draw = none},tick style={draw = none}]
\pgfmathsetmacro{\a}{0}
\pgfmathsetmacro{\b}{pi/3}

\addplot3[surf,
      samples=20,
      domain=0:2*pi,y domain=0:pi,
      z buffer=sort,
      opacity = 0.4,
      fill = black!40!white,
      colorbar,
      colormap ={grayscale}{gray(0cm)=(0.5); gray(1cm)=(0.5)},
      ]
      ({spherepar1(x,y)},{spherepar2(x,y)},{spherepar3(x,y)});
     \addplot3[line width = 2pt,
    opacity = 0.6,
      samples=50,
      domain=0:2*pi,
      samples y = 0,
      ]
      ({sin(deg(x))},{0},{cos(deg(x))});
\addplot3[surf,
      samples=10,
      domain=-1:1,y domain=-1:1,
      z buffer=sort,
      opacity = 0.4,
      fill = blue!40!white,
      colorbar,
      colormap ={grayscale}{gray(0cm)=(0.5); gray(1cm)=(0.5)},
      ]
      ({spherepar1(\a,\b) +x*der11(\a,\b) +y*der21(\a,\b)}, 
      {spherepar2(\a,\b) +x*der12(\a,\b) +y*der22(\a,\b)},
        {spherepar3(\a,\b) +x*der13(\a,\b) +y*der23(\a,\b)}
        );
        \node at (axis cs:{spherepar1(\a,\b)}, {spherepar2(\a,\b)},{spherepar3(\a,\b)}) [circle, scale=0.3, draw=black!100,fill=black!100] {};
        \node at (axis cs:{spherepar1(\a,\b)}, {spherepar2(\a,\b)},{spherepar3(\a,\b)}) [above] {$x_0$};
\addplot3[line width = 2pt,
      samples=10,
        samples y =0,
      domain=-1:1,
      z buffer=sort,
      opacity = 0.6,
      color = blue!80!white,
      ]
      ({spherepar1(\a,\b) +x*cos(deg(\b)) }, 
      {spherepar2(\a,\b) },
        {spherepar3(\a,\b) -x*sin(deg(\b))}
        );
        \node at (axis cs:0.6,0.75,1.1)  {$T_{x_0}N \subset  T_{x_0}M$};
        \node at (axis cs:-1,0,0.7)  {$N$};
\end{axis}
\end{tikzpicture}
\begin{tikzpicture}[scale = 1.2,
declare function = {spherepar1(\x,\y) = 
{sin(deg(\y))*cos(deg(\x))}; 
},
declare function = {
spherepar2(\x,\y) = {sin(deg(\y))*sin(deg(\x))};
},
declare function = {
spherepar3(\x,\y) ={cos(deg(\y))};
},
declare function = {
der11(\x,\y) ={sin(deg(\x))};
},
declare function = {
der12(\x,\y) ={cos(deg(\x))};
},
declare function = {
der13(\x,\y) ={0};
},
declare function = {
der21(\x,\y) ={cos(deg(\x))*cos(deg(\y))*(1-2*sin(deg(\x))^2)/sqrt(1 + 4*cos(deg(\y))^2*cos(deg(\x))^4 -4*cos(deg(\y))^2*cos(deg(\x))^2)};
},
declare function = {
der22(\x,\y) ={sin(deg(\x))*cos(deg(\y))*(1-2*cos(deg(\x))^2)/sqrt(1 + 4*cos(deg(\y))^2*cos(deg(\x))^4 -4*cos(deg(\y))^2*cos(deg(\x))^2)};
},
declare function = {
der23(\x,\y) ={-sin(deg(\y))/sqrt(1 + 4*cos(deg(\y))^2*cos(deg(\x))^4 -4*cos(deg(\y))^2*cos(deg(\x))^2)};
},
]
\tikzstyle{every node}=[font=\small]
\begin{axis}[axis equal image,ymin = -1.25,ymax = 1.25,zmin = -1.25,zmax = 1.25,xmin=-1.25,xmax=1.25,axis line style={draw = none},tick style={draw = none},ticks = none,colorbar,
      colorbar style = {ylabel = f(x),height = 3cm,width = 0.2cm,at = {(1.03,0.5)},anchor = west}]
\pgfmathsetmacro{\a}{0.}
\pgfmathsetmacro{\b}{pi/3}
\addplot3[surf,
      samples=50,
      domain=0:2*pi,y domain=0:pi,
      z buffer=sort,
      opacity = 0.4,
      point meta = {sin(deg(y))^2},
      ]
      ({spherepar1(x,y)},{spherepar2(x,y)},{spherepar3(x,y)});
\addplot3[line width = 2pt,
    opacity = 0.4,
      samples=50,
      domain=0:2*pi,
      samples y = 0,
      ]
      ({sin(deg(x))},{0},{cos(deg(x))});
        \node at (axis cs:{spherepar1(\a,\b)}, {spherepar2(\a,\b)},{spherepar3(\a,\b)}) [circle, scale=0.3, draw=black!100,fill=black!100] {};
        \node at (axis cs:{spherepar1(\a,\b)}, {spherepar2(\a,\b)},{spherepar3(\a,\b)}) [above] {$x_0$};
        \node at (axis cs:-1,0,0.7)  {$\Z$};
\end{axis}
\end{tikzpicture}

\vspace*{-.75cm}

\caption{\textbf{Left:} Representation of the manifold $M = S^2$ as well as a sub-manifold $N$ homeomorphic to a circle. The tangent spaces at a given point $x_0 \in N\subset M$ are represented as well. \textbf{Right:} Representation of a non-negative function on the sphere as a color map; it satisfies the NHC, and its null space $\Z$ is represented in black.  }
\end{figure}

For more precise definitions of topological spaces, atlases, charts, and details on the ``Hausdorff second-countable'' condition, see for instance \cite{lafontaine2015introduction,paulingeodiff,topologie,spivak}. The main idea behind the introduction of manifolds as opposed to sub-manifolds of $\R^d$ is to consider the intrinsic geometric object, and not its relation to the euclidean space it is embedded in (as such an embedding is not unique). An example of manifold as well as a representation of the tangent space is provided in the left of \cref{fig:manifold_fig}. 

\subsection{Assumptions in the manifold case}\label{sec:reformulation_assumptions}
In this section, we formulate the local NHC in the case of manifolds, and rewrite \cref{thm:local_decomposition} in this setting. We also extend the definitions of being positively normal and the global NHC (\cref{df:global_nhc}).

Fix $p \in \N,~p \geq 2$, and a manifold $M$ of regularity at least $C^p$ and of dimension $d \in \N$. 
To start with, let $f: \Omega \rightarrow \R$ be a non-negative function defined on an open set of $M$ and of class~$C^p$. As before, define $\Z$ to be the set of zeros of $f$.

Contrary to the $\R^d$ case, the second differential of the function $f$ cannot be identified to a symmetric bi-linear form everywhere. However, it is the case at so-called critical points, i.e.,  points $x \in \Omega$ such that $d f(x) = 0$. In particular, since all the zeros of a $C^1$ non-negative function are critical points, this Hessian will be defined at all points in the set of zeros $\Z$ of~$f$.

\begin{lemma}[definition of the Hessian]
Let $x$ be a critical point of $f$. Then there exists a unique symmetric bi-linear form $H_f(x): T_x M \times T_x M \rightarrow \R$ such that for any local chart $\phi: (M,x) \rightarrow (\R^d,0)$ it holds:
\[\forall \xi,\eta \in T_x M \times T_x M,~ H_f(x)[\xi,\eta] = d^2 (f\circ \phi^{-1})(0)[d\phi(x)\xi,d\phi(x)\eta]  . \]
\end{lemma}
In order to prove this lemma, we simply define the bilinear form as such for a given chart $\phi$ around $x$, and then show that this definition does not depend on the chart $\phi$ using the fact that $x$ is a critical point. This is completely proved in section 2 of \cite{milnor1963}. 
In order to formulate the definition of the normal Hessian condition in the setting of manifolds, we further need a definition of what a sub-manifold of $M$ is.
A subset $N \subset M$ is said to be a sub-manifold of $M$ of class $C^k$ if $M$ is of class $C^k$ and if, for any $x \in N$ and any local chart $\phi: U \rightarrow \R^d$ defined on a neighborhood of $x$, $\phi(U\cap N)$ is a sub-manifold of $\R^d$ of class $C^k$. In the literature, this is also called a proper sub-manifold (see on the left of \cref{fig:manifold_fig} for an example).

\begin{definition}[Normal Hessian condition for a manifold]\label{df:NHC_manifold} Let $x \in \Omega$ be a point in the domain of $f$. We say that $f$ satisfies the normal Hessian condition (NHC) at $x$ if there exists a dimension $0 \leq d_0 \leq d$ and a sub-manifold $N$ of $M$ of class $C^k$ with $k \geq 1$ and of dimension $d_0$ such that $x \in M \subset \Z$, and one of the following equivalent conditions is satisfied: 

\begin{itemize}
    \item[(i)] the rank of $H_{f}(x)$ is $d-d_0$;
    \item[(ii)] $H_{f}(x)[\xi,\xi] >0$ for any vector $\xi \in T_x M \setminus{T_x N}$.
\end{itemize}
\end{definition}

Using any local chart around the point $x$ in the domain of $f$, one can apply \cref{thm:local_decomposition} to obtain the following theorem as a corollary.

\begin{theorem}\label{thm:local_decomposition_manifold}
If $f$ satisfies the NHC at $x_0$ with regularity $k$ and dimension $d_0$, there exists an open neighborhood $U$ of $x_0$ in $M$ on which $f$ is defined and such that $U \cap \Z$ is a sub-manifold of $M$ of dimension $d_0$ and of class $C^{\max(k,p-1)}$; and there exist functions $f_i \in C^{p-2}(U)$ where $1 \leq i \leq d-d_0$ such that 
\begin{equation}
    \label{eq:local_form_manifold}
    \forall x \in U,~ f(x) = \sum_{i=1}^{d-d_0}{f_i^2(x)}.
\end{equation}
\end{theorem}

Exactly in the same way as for the definition of the NHC for manifolds, we can similarly extend the definition of a function being positively normal to a sub-manifold  in \cref{df:positive_normal} as well as the global NHC in \cref{df:global_nhc}. We will therefore say that $f: M \rightarrow \R$ which is non-negative satisfies the global NHC if it satisfies the local NHC at every point in its set of zeros $\Z$, or equivalently if it is positive normally to $\Z$ which is a sub-manifold of $M$ of class $C^1$ (or $C^{p-1}$). 
On the right hand side of \cref{fig:manifold_fig}, we represent a function which satisfies the NHC on the sphere $S^2$ through a colormap, with a continuous set of zeros. The goal is to prove that such a function can be decomposed as a sum of squares on $S^2$.

\subsection{Gluing local decompositions to form a global one}\label{sec:gluing}

In this section, we present and develop the tools to glue local decompositions such as \cref{thm:local_decomposition_manifold} into a global one, which will lead to \cref{thm:main_manifold}.

The first result we need is a simple result to ``extend'' a function defined on an open set $U$ of~$M$ to $M$ by multiplying it by a function defined on $M$ whose support lies in $U$ (\cref{lm:extension_lemma}). The second one is a variant of the fundamental result of existence of partitions of unity on a manifold, adapted to our sum of squares setting (\cref{lm:glue_square}). Recall that the support of a function has been defined in \cref{secsec:global_sos_rd}. The proof of these results can be found in \cref{app:paracompactness}.

\begin{lemma}[Extension lemma]\label{lm:extension_lemma}
Let $q \in \N$, $M$ be a manifold of class at least $C^q$. Let $U$ be an open set of $M$, $g : U \rightarrow \R$ be a $C^q$ function defined on $U$, and $\chi : M \rightarrow \R$ be a $C^q$ function defined on the whole of $M$ but with support included in $U$. 
Then the function $\chi g : U \rightarrow \R$ extended as $0$ on $M\setminus{U}$, is of class $C^q$ on the whole of $M$ and has support included in $\supp(\chi) \subset U$. We still denote with $\chi g$ its extension to $M$.
\end{lemma}

\begin{lemma}[Gluing lemma]\label{lm:glue_square}
Let $(U_i)_{i \in I}$ be an open covering of a manifold $M$ of class $C^k$ (i.e. $\bigcup_{i \in I}{U_i} = M$). There exists a family of functions $\chi_i: M \rightarrow [0,1]$ of class $C^k$ with locally finite support, such that $\supp(\chi_i) \subset U_i$ for all $i \in I$ and satisfying: 
\[\sum_{i \in I}{\chi_i^2} = 1.\]
\end{lemma}


We can now proceed from local to global in two steps. First, we use the gluing lemma to glue decompositions in a single connected component of the manifold of zeros (\cref{lm:glue_connected}). We then glue these different decompositions into a single global one (\cref{lm:glue_connected_together}).  

\begin{lemma}\label{lm:glue_connected} Assume $f$ satisfies the global NHC. Let $N$ be a connected component of its manifold of zeros $\Z$. There exists an open neighborhood $U$ of $N$ as well as a locally finite, at most countable family $(f_j)_{j \in J}$ of functions of class $C^{p-2}$ such that
\begin{equation}\label{eq:sos_one_connected}
\forall x \in U,~ f(x) = \sum_{j \in J}{f_j(x)^2}.
\end{equation}
Moreover, we can find $J$ such that \textit{a)} $|J| = d$ if $N= \{x_0\}$ is a single point and \textit{b)} $J$ is finite if $N$ is compact.

\end{lemma}

\begin{proof}
The case where $N = \{x_0\}$ is simply \cref{thm:local_decomposition_manifold} applied to $x_0$. In the other cases, note that for all $x \in N$, by \cref{thm:local_decomposition} since the NHC is satisfied at $x$, there exists an open neighborhood $U_x$ of $x$ as well as functions $(f_{x,i})_{1 \leq i \leq d}$ of class $C^{p-2}$ such that $f = \sum_{i=1}^d{f_{x,i}^2}$ on $U_x$. Since   $(U_x)_{x \in N}$ covers $N$, we can extract a covering $(U_{x_j})_{j \in J}$ of $N$ such that \textit{a)} $J$ is finite if $N$ is compact and \textit{b)} $J$ is at most countable otherwise, since $N$ is second-countable and Hausdorff. Denote with $(U_{j})_{j \in J}$ this open covering, and replace $x$ by $j$ to denote the associated $f_{j,i}$. Denote with $U$ the open set $\bigcup_{j}{U_j}$.

Applying \cref{lm:glue_square} to the manifold $U$, we can find a family of functions $(\chi_j)_{j\in J}$ with locally finite support, such that $\supp(\chi_j) \subset U_j$ and $\sum_{j}{\chi_j^2} =1$ on $U$. By the extension \cref{lm:extension_lemma}, we can therefore define the functions $\widetilde{f}_{j,i} := \chi_j~f_{j,i}$ for $i \in \{1,...,d\}$ and $j \in J$ which are defined on the whole of $M$. Note that since $\supp(\tilde{f}_{j,i}) \subset \supp(\chi_j)$ and since $1 \leq i \leq d$, the support of $(\tilde{f}_{j,i})$ is also locally finite. To conclude, we use the property that $\sum_{j}{\chi_j^2 } = 1$ on $U$ as well as the fact that $\sum_{i}{f_{j,i}^2} = f$ on $\supp(\chi_j) \subset U_j$ to show that $\sum_{i,j}{\tilde{f}_{j,i}^2} = f$ on $U$. The number of functions $\tilde{f}_{j,i}$ is finite if $N$ is compact since $J$ is finite, and is at most countable else since $J$ is at most countable.
\end{proof}

\begin{lemma}\label{lm:glue_connected_together}
Let $\Z = \sqcup_{i \in I}{N_i}$ be the manifold of zeros decomposed along its connected components. 
Assume that there exists an index set $J$, such that for all $i \in I$, there exists an open neighborhood $U_i$ of $N_i$ on which $f$ can be decomposed as a sum of squares indexed by $J$:
\begin{equation}\label{asm:lm:connected}
\forall i \in I,~ \exists (f_{i,j})_{ j \in J} \in (C^{p-2}(U_i))^J,~\forall x \in U_i,~ f(x)  = \sum_{j \in J}{f_{i,j}(x)^2},
\end{equation}
and such that the families $(f_{i,j})_{ j \in J}$ are all locally finite. 
Then there exists a locally finite family $(g_j)_{j \in J\cup \{\star\}}$ of $C^{p-2}$ functions on $M$ (we add an extra element $\star$ to $J$), such that 
\begin{equation}\label{eq:res_lemma_connected}
\forall x \in M,~f(x) = \sum_{j \in J\cup\{\star\}}{g_j(x)^2}.
\end{equation}
\end{lemma}

\begin{proof}
By \cref{lm:open_sets}, there exist disjoint open sets $V_i\subset M$ such that $N_i \subset V_i$, since $\Z$ is a proper sub-manifold of $M$ by \cref{df:global_nhc} (directly adapted to the manifold case). Hence, we can assume that the $U_i$ are disjoint (consider instead $U_i \cap V_i$, the property still holds). 
Define $U_{\star} = \{f > 0\}$. Note that since the $U_i$'s cover $\Z$, $U_{\star} \cup \bigcup_{i \in I}{U_i}$ covers $M$ since $f$ is non-negative: take $\chi_{\star},(\chi_i)_{i \in I}$ to be a gluing family adapted to that covering given by \cref{lm:glue_square}.
Note that for any $i,i^{\prime} \in I$, we have $\chi_i \chi_{i^{\prime}} = 0$ since $\chi_i$ is supported on $U_i$ and the $U_i$ are disjoint. 
Consider the function $g_j = \sum_{i \in I}{\chi_i f_{i,j}}$, which is well defined on $M$ and $C^{p-2}$ by \cref{lm:extension_lemma}. We have $g_j^2 = \sum_{i \in I}{\chi_i^2 f_{i,j}^2}$ since $\chi_i \chi_{i^{\prime}} = 0$ when $i \neq i^{\prime}$. 

\noindent \textit{Assertion : the family $(g_j)_{j\in J}$ has locally finite support.} 
Let $x \in \R^d$ and assume $g_j(x) \neq 0$. Then there exists $i \in I$ such that $\chi_i(x) >0$, and hence there exists an open set $U_x$ around $x$ such that $U_x \subset U_i$. But in that case, $\chi_{i^{\prime}}(x^{\prime}) = 0$ for all other $i^{\prime}$ and for all $x^{\prime} \in U_x$ since the $U_i$ are disjoint and $\chi_{i^{\prime}}$ is supported on $U_{i^{\prime}}$. Moreover, since $(f_{i,j})_{j \in J}$ is locally finite, there exists an open set $V_x$ around $x$ as well as a finite $J_0\subset J$ such that $f_{i,j} = 0$ on $V_x$ for all $j \in J\setminus{J_0}$. Hence, for any $j \in J\setminus{J_0}$ and any $x^{\prime} \in U_x \cap V_x$, we have $f_{i,j}(x^{\prime})= 0$ and $\chi_{i^{\prime}}(x^{\prime}) = 0$ thus $g_j(x^{\prime}) = 0$. Thus, $U_x\cap V_x \subset M\setminus{\supp(g_j)}$ for all $j \notin J_0$: the family $g_j$ is locally finite.

\noindent \textit{Conclusion.} Define $g_{\star} = \chi_{\star} \sqrt{f}$, which is of class $C^{p}$ since $\chi_{\star}$ is supported on $\{f>0\}$. Since the addition of one function changes nothing to the locally finite property of a family of functions, the family $(g_j)_{j \in J}\cup g_{\star} $ is still locally finite. Using the fact that $g_j^2 = \sum_{i \in I}{\chi_i^2 f_{i,j}^2}$, that $\sum_{j \in I\cup \{\star\}}{\chi_i^2} =1$ and \cref{asm:lm:connected}, it holds
\begin{align*}
    \sum_{j \in J\cup \{\star\}}{g_j^2} &= \chi^2_{\star} f + \sum_{j \in J}{g_j^2}
     = \chi^2_{\star} f + \sum_{j \in J}{\sum_{i \in I}{\chi_i^2f_{i,j}^2}}\\
    &= \chi^2_{\star} f + \sum_{i \in I}{\sum_{j \in J}{\chi_i^2f_{i,j}^2}}
    = \chi^2_{\star} f + \sum_{i \in I}{\chi_i^2 f} = f.
\end{align*}

\end{proof}

\subsection{Main results}\label{sec:general_manifold}
We are now ready to state our main result on manifolds. On the right hand side of \cref{fig:manifold_fig}, we represent a case where this theorem applies for a non-negative function defined on $S^2$.

\begin{theorem}
\label{thm:main_manifold} 
Let $M$ be a manifold and $f: M \rightarrow \R$ be a non-negative map of class $C^p$.
Assume $f$ satisfies the global normal Hessian condition. Then there exists $I$ which is at most countable and functions $f_i \in C^{p-2}(M)$ for $i \in I$ such that the family $(f_i)$ has locally finite support and  
\begin{equation}
    \label{eq:res_manifold}
    \forall x \in M,~ f(x) = \sum_{i\in I}{f_i(x)^2}.
\end{equation}
Moreover: 
\begin{itemize}
    \item if $f$ satisfies the strict Hessian condition, $\Z$ is discrete and we can find such a decomposition such that $|I| \leq d+1$.
    \item if $\Z$ is compact, then $|I|$ can be taken to be finite.
    \end{itemize}
\end{theorem}

\begin{proof}
The proof of this theorem is a simple consequence of \cref{lm:glue_connected} and \cref{lm:glue_connected_together}. Note that the global NHC \cref{df:global_nhc} shows that $\Z$ is a sub-manifold of $M$. Let $N_i$ denote the connected components of $\Z$. By \cref{lm:open_sets}, we can find disjoints open sets $U_i$ such that $N_i \subset U_i$. 

\noindent \textit{General case.} Without any more assumptions, we know from \cref{lm:glue_connected} that on any connected component $N_i$, we can have a decomposition of the form $f = \sum_{j \in J}{f_{i,j}^2}$ with $f_{i,j}\in C^{p-2}$ a family with locally finite support on an open neighborhood $V_i$ of $N_i$. Moreover, we know that $J$ is at most countable. Adding zeros when necessary, and reindexing, we can assume that $J = \N$. Now applying \cref{lm:glue_connected_together}, we prove the general case. 

\noindent \textit{Compact case.} If we assume that $N$ is compact, since the $U_i$ cover $N$, necessarily the number of connected components is finite (just extract a finite covering of $N$ from the $U_i$). We know from \cref{lm:glue_connected} that on any connected component $N_i$, we can have a decomposition of the form $f = \sum_{j = 1}^{n_i}{f_{i,j}^2}$ with $f_{i,j}\in C^{p-2}$ and $n_i \in \N$ on an open neighborhood $V_i$ of $N_i$, since $N_i$ is compact. Hence, up to adding $f_{i,j} = 0$, we can assume that $n_i = n = \max_{i}(n_i)$ since there are a finite number of connected components. Now applying \cref{lm:glue_connected_together} with $J = \{1,...,n\}$, the result is proven in the compact case with $n+1$ functions.

\noindent \textit{SHC case.} If we assume that the SHC holds, every connected component $N_i$ is a singleton $\{x_i\}$: we know from \cref{lm:glue_connected} we can have a decomposition of the form $f = \sum_{j = 1}^{d}{f_{i,j}^2}$ with $f_{i,j}\in C^{p-2}$ on an open neighborhood $V_i$ of $N_i$, since $N_i$ is compact.  Now applying \cref{lm:glue_connected_together} with $J = \{1,...,d\}$, the result is proven with $d+1$ functions.

\end{proof}

\begin{remark}
Note that the difference between the number of functions in the SHC case is better than the one obtained in \cite{rudi2020finding}. This is because of the two step procedure in the gluing: first in a connected component, and then between connected components. The long term goal is to be able to prove that we need only a finite number $N(d)$ of functions per connected component (in the compact case), and hence to have an explicit bound after gluing the connected components together, rather than just relying on a compact extraction argument, which is not as precise.
\end{remark}

\section{Proof of the local decomposition as a sum of squares}\label{sec:proofs}\label{sec:proof_local_decomposition}

In this section, we formally prove the key result of the paper, \cref{thm:local_decomposition}.

\begin{proof}

Note that the existence of the sub-manifold $N$ of dimension $d_0$ around $x_0$ implies the existence of a local parametrization around $x_0$ (see for instance Theorem 1.21 of \cite{lafontaine2015introduction}): there exists an open neighborhood $\Wt_0$ of $0$ in $\R^{d_0}$, an open neighborhood $U_{x_0}$ of $x_0$ in $\R^d$ and  a $C^k$ immersion $\phi: \Wt_0 \rightarrow U_{x_0}$ of class $C^k$ such that $\phi$ is a homeomorphism from $\Wt_0$ onto $U_{x_0} \cap N$. Since restricting $N$ to $N \cap U_{x_0}$ does not change the assumptions of the theorem, we will will assume that  $N = \im(\phi)$ for a $C^k$ immersion $\phi: (0,\R^{d_0}) \rightarrow (x_0,\R^d)$. We will denote with $T_{x_0}:= d\phi_0(\R^{d_0})= T_{x_0}N$ the tangent space to $N$ at $x_0$.

Before starting the proof, recall that for any $x \in \Z$, it holds $df(x) = 0$ and $d^2 f (x) \succeq 0$ (or equivalently $\nabla^2 f(x) \succeq 0$). Moreover, note that if $A \in \psd(\R^d)$ is a symmetric positive semi-definite matrix, if a vector $k \in \R^d$ satisfies $k^{\top} A k = 0$, then $A k = 0$ (this is a trivial consequence of the spectral theorem by decomposing $k$ along an orthonormal basis of eigenvectors).

\noindent \textit{Step 1: characterizing the null-space of the Hessian.}
We will prove that under the assumptions of the theorem, \textit{a)} $T_{x_0}$ is equal to the null-space $\ker( \nabla^2 f(x_0))$ of the Hessian of $f$ at $x_0$ and \textit{b)} that for any supplementary $S$ to $T_{x_0}$, the restricted Hessian $\nabla^2 f(x_0)|_{S}$ is positive definite.

To prove \textit{a)}, assume that there exists an element in $ k \in T_{x_0}$ such that $\nabla^2f(x_0) k \neq 0$. Since $\nabla^2 f(x_0)$ is positive semi-definite, this implies that $k^{\top} \nabla^2 f(x_0) k > 0$. Let $h \in \R^{d_0}$ such that $d \phi_0 h = k$, and let $x_t = \phi(th)$ which is defined for $t$ in an open neighborhood of $0$. Using the Taylor expansion of $f$ around $x_0$:
\[f(x) - f(x_0) - d f(x_0)[x-x_0]= \tfrac{1}{2}(x-x_0)^{\top}\nabla^2 f(x_0) (x-x_0) + \epsilon(x-x_0)\|x-x_0\|^2,\]
where $\epsilon(x) \underset{\|x\|\rightarrow 0}{\rightarrow} 0$. Now applying this for $x_t$, since $f(x_t) = f(x_0) = 0 $ and $df(x_0) = 0$, it holds:
\[0 = \tfrac{1}{2}(x_t-x_0)^{\top}\nabla^2 f(x_0) (x_t-x_0) + \epsilon(x_t-x_0)\|x_t-x_0\|^2.\]
Using the fact that $\phi$ is differentiable at $0$ yields $x_t - x_0 = t d\phi(0)[h] + o_{t \rightarrow 0}(t) = tk + o_{t \rightarrow 0}(t)$. Injecting this in the equation above yields
\[0 = t^2 ~ \tfrac{1}{2}k^{\top}\nabla^2 f(x_0) k + o_{t \rightarrow 0}(t^2).\]
Hence, necessarily, $k^{\top}\nabla^2 f(x_0) k = 0$, which is a contradiction. This proves that $T_{x_0} \subset \ker(\nabla^2 f (x_0))$, and in particular, $d_0 \leq \dim(\ker(\nabla^2 f (x_0)))$. Since the rank of $\nabla^2 f(x_0)$ is actually $d-d_0$, the rank theorem shows that $\dim(\ker(\nabla^2 f (x_0))) = d_0$ and hence $T_{x_0} = \ker(\nabla^2 f (x_0))$.

To prove \textit{b)}, we just need to prove that the restriction to any supplementary to the null-space of $\nabla^2 f(x_0)$ is positive definite. Using the small result at the beginning of the proof, any vector $k \in \R^d \setminus{T_{x_0}}$ satisfies $k^{\top} \nabla^2 f(x_0) k > 0$. In particular, this means that the restriction of $\nabla^2 f(x_0)$ to any supplementary subspace $S$ of $T_{x_0}$ is positive definite.

\noindent \textit{Step 2: applying the Morse lemma. }
Let $P = (P_1,P_2) \in O_d(\R)$ be the matrix of an orthonormal basis adapted to the decomposition $\R^d =  T_{x_0}^{\perp} \oplus T_{x_0}$. Note that $P_1 \in \R^{d \times (d-d_0)}$ and $P_2 \in \R^{d \times d_0}$ are also orthonormal matrices, and that since $P_1$ spans $T_{x_0}^{\perp}$, in particular $P_1^{\top}\nabla^2 f(x_0)P_1 \succ 0$. 

Define $g: (\xp,\yp) \in \R^{d-d_0} \times \R^{d_0} \mapsto  f(P_1\xp  + P_2 \yp + x_0) = f(\Aa(\xp,\yp))$, where $\Aa(\xp,\yp) = P (\xp,\yp) + x_0$ is an isometry\footnote{An isometry is simply a map which preserves distances, and can be defined as an orthogonal transformation plus an affine shift.} ($\Aa^{-1}x = P^{\top}(x-x_0)$). We have $\nabla_{\xp} g(0,0) = P_1^{\top} \nabla f(x_0) = 0$ and $\nabla^2_{\xp\xp} g(0,0) = P_1^{\top} \nabla^2 f(x_0) P_1 \succ 0$. We can therefore apply the Morse lemma \cref{lm:morse_hormander} to $g$: there exists two open neighborhoods of zero $V \subset \R^{d-d_{0}},W\subset \R^{d_0}$ as well as $\varphi: W \rightarrow V$ of class $C^{p-1}$ such that $\{(\xp,\yp) \in V \times W~:~\nabla_{\xp} g(\xp,\yp) = 0\} = \{(\xp,\yp) \in V \times W~:~ \xp = \varphi(\yp)\}$ and $z: V \times W \rightarrow \R^{d-d_0}$ of class $C^{p-2}$ such that
\begin{equation}
    \label{eq:result_morse}
    \forall (\xp,\yp) \in V \times W,~ g(\xp,\yp) = g(\varphi(\yp),\yp) + \tfrac{1}{2}z(\xp,\yp)^{\top} H^{\prime} z(\xp,\yp),
\end{equation}
where $H^{\prime}$ is the positive definite matrix  $P_1^{\top} \nabla^2 f(x_0) P_1$. 

\noindent\textit{Step 3: making the sum of squares appear.}
Since $H^{\prime} \in \psd(\R^{d-d_0})$ and $H^{\prime} \succ 0$, we can decompose it using the spectral theorem: $H^{\prime} = \sum_{i=1}^{d-d_0}{\lambda_i u_i u_i^{\top}}$ where the $\lambda_i > 0$. Defining $g_i(\xp,\yp) = \sqrt{\lambda_i/2} u_i^{\top}z(\xp,\yp) $, \cref{eq:result_morse} can be rewritten as 
\begin{equation}
    \label{eq:result_morse_bis}
    \forall (\xp,\yp) \in V \times W,~ g(\xp,\yp) = g(\varphi(\yp),\yp) + \sum_{i=1}^{d-d_0}{g_i^2(\xp,\yp)}.
\end{equation}
Note that the $g_i$ are of class $C^{p-2}$ since $z$ is of class $C^{p-2}$. 
We see that if we can show that $g(\varphi(\yp),\yp)=0$ in a neighborhood of $(0,0)$, since we can go back to the original coordinate system through $\Aa^{-1}$, we will have shown the theorem.

\textit{Step 4: characterizing $\Z$ in a neighborhood of $x_0$.} 
Denote with $G_{\varphi} = \{(\varphi(y),y)~:~ y \in W\}$ the graph of $\varphi$, and which is a sub-manifold of class $C^{p-1}$ of $\R^{d-d_0} \times \R^d$ (see theorem 1.21, point (iv) of \cite{lafontaine2015introduction}). Since $\Aa$ is an isometry, the set $\Aa(G_{\varphi})$ is also a sub-manifold of class $C^{p-1}$ of $\R^d$.

Let $\Wt = \phi^{-1}(\Aa(V \times W))$: it is an open neighborhood of $0$.
Note that $\phi(\Wt) \subset \Z \cap \Aa(V \times W)$ by assumption, and since for any $x \in \Z$, we have $\nabla f(x) = 0$, it holds in particular that for any $x \in \Z \cap \Aa(V\times W)$, we have $\nabla_{\xp}g(\Aa^{-1}(x)) = P_1^{\top} \nabla f(x) = 0$. Hence, by the result of the Morse lemma, it holds $\Aa^{-1}(\phi(\Wt)) \subset \Aa^{-1}(\Z)\cap (V \times W) \subset G_{\varphi}$.

Define $\psi: (\xp,\yp) \in V\times W \mapsto (\xp-\varphi(\yp),\yp)$ which is a $C^{p-1}$ diffeomorphism onto its image with inverse $(t,u) \mapsto (t + \varphi(u),u)$. Note that $\psi$ maps $G_{\varphi}$ onto $\{0_{\R^{d-d_0}}\} \times W$. If $\pi_2$ denotes the canonical projection $\pi_2:\R^{d-d_0} \times \R^{d_0} \rightarrow \R^{d_0}$, we see that $\pi_2 \circ \psi$ maps $G_{\varphi}$ injectively onto $W \subset \R^{d_0}$. 

 Take $\Phi = \pi_2 \circ \psi \circ \Aa^{-1} \circ \phi: \widetilde{W} \rightarrow \R^{d_0}$, which is well defined by definition of $\Wt$, and $C^1$ by composition. Note that it is an immersion at $0$. Indeed \textit{i)} $\phi$ maps $0$ onto $x_0$ and is an immersion at $0$ by assumption, hence $d \phi_0$ is injective; \textit{ii)} $\psi \circ \Aa^{-1}$ is a $C^{p-1}$ diffeomorphism from $\Aa(V\times W)$ (containing $x_0$) to its image, and hence its differential is invertible at $x_0$, and thus by composition, the  differential $d (\psi \circ \Aa^{-1} \circ \phi)(0)$ is injective; \textit{iii)} since $\Aa^{-1}(\phi(\Wt)) \subset G_{\varphi}$ by a previous statement, and since $\psi(G_{\varphi}) \subset \{0\} \times W$ also by a previous statement, it, it holds that the differential $d (\psi \circ \Aa^{-1} \circ \phi)(0) \R^{d_0} \subset \{0\}\times \R^{d_0}$ and hence applying $\pi_2$ does not change the injectivity of the differential; hence $\Phi$ is an immersion at $0$. But since $d \Phi_0$ is a linear map from $\R^{d_0}$ to $\R^{d_0}$, $d\Phi_0$ being injective is equivalent to $d \Phi_0$ being invertible. Hence, by the local inversion theorem \cref{thm:loc_inversion}, there exists an open neighborhood of $0$ $\widetilde{W}^{\prime} \subset \Wt$ and an open neighborhood of $0$ $W^{\prime} \subset W$ such that $\Phi$ is a $C^1$ diffeomorphism from $\Wt^{\prime}$ to $W^{\prime}$. 
 
 Define $U = (\pi_2 \circ \psi \circ  \Aa^{-1})^{-1}(W^{\prime}) = \Aa( \psi^{-1}(\R^{d-d_0}\times W^{\prime}))$ ,which is an open neighborhood of $x_0$. Note that since $\Phi$ is a diffeomorphism from $\Wt^{\prime}$ to $W^{\prime}$, we have $\phi(\Wt^{\prime}) \subset U$. Moreover, since $\psi$ is defined on $V\times W$, we have $U \subset \Aa(V\times W)$. Finally, let $u \in U \cap \Aa(G_{\varphi})$. Since $u \in U$, there exists $\tilde{w}^{\prime} \in \Wt^{\prime}$ such that $\pi_2\circ \psi \circ \Aa^{-1}(\phi(\tilde{w}^{\prime})) = \pi_2\circ \psi \circ \Aa^{-1}(u)$. Moreover, since $\pi_2 \circ \psi$ is injective on $G_{\varphi}$, and since both $\Aa^{-1}(\phi(\tilde{w}^{\prime}))$ and $\Aa^{-1}(u)$ belong to $G_{\varphi}$ (the first using the previous point since $\Wt^{\prime}\subset \Wt$ and the second by assumption), we have  $\Aa^{-1}(\phi(\tilde{w}^{\prime}))=\Aa^{-1}(u)$ and hence $u  = \phi(\tilde{w}^{\prime})$ since $\Aa$ is one to one. This shows that $U \cap \Aa(G_{\varphi}) \subset \phi(\Wt^{\prime})$.
 
Moreover, a previous point shows that $\Aa^{-1}(\phi(\Wt)) \subset \Aa^{-1}(\Z)\cap (V \times W) \subset G_{\varphi}$. Now since $\Aa$ is one to one and since $\Wt^{\prime} \subset \Wt$ we have $ \phi(\Wt)\subset  \Z \cap (\Aa(V \times W)) \subset \Aa(G_{\varphi})$. Since $\phi(\Wt^{\prime}) \subset U$, we therefore have $\phi(\Wt^{\prime}) \subset \Z \cap U \subset \Aa(G_{\varphi}) \cap U$. Combining this with the previous result, we finally have 
\begin{equation}\label{eq:equality}
\phi(\Wt^{\prime}) \subset \Z \cap U \subset \Aa(G_{\varphi}) \cap U \subset \phi(\Wt^{\prime}) \implies \phi(\Wt^{\prime}) = \Z \cap U =\Aa(G_{\varphi}) \cap U.
\end{equation}
\noindent \textit{Step 5: conclusion.}
\cref{eq:equality} shows that $\phi(\Wt^{\prime}) =\Z \cap U = \Aa(G_{\varphi}) \cap U$.

One the one hand, this shows that $U \cap \Z$ is the intersection between an open set $U$ and a  sub-manifold $\Aa(G_{\varphi})$ of $\R^d$ of class $C^{p-1}$ (since it is the composition of the graph of $\varphi$ which is $C^{p-1}$, which is a $C^{p-1}$ manifold by \cite{lafontaine2015introduction}, by an isometry which is in particular a diffeomorphism). Moreover, since $\phi$ is a $C^k$ immersion which is a homeomorphism on its image, $\phi(\Wt^{\prime})$ is a sub-manifold of class $C^k$. Thus, $U\cap \Z$ is a sub-manifold of $\R^d$ of class $C^{\max(k,p-1)}$.

On the other, since $\Aa^{-1}(U) \subset V \times W$, \cref{eq:result_morse_bis} becomes 
\begin{equation}\label{eq:intermediate}
    \forall u \in U,~ g(\Aa^{-1}(u)) = g(\varphi(\yp),\yp) + \sum_{i=1}^{d-d_0}{g_i^2(\Aa^{-1}(u))},~(\xp,\yp) = \Aa^{-1}(u).
\end{equation}
Let $u \in U$ and write $(\xp,\yp) = \Aa^{-1}(u)$. First, note that $\Aa(\varphi(\yp),\yp) \in \Aa(G_{\varphi})$. Moreover, since $\Aa^{-1}u \in \psi^{-1}(\R^{d-d_0}\times W^{\prime})$ by definition of $U$, this shows that $\yp \in W^{\prime}$ and hence $(\varphi(\yp),\yp) = \psi^{-1}(0,\yp) \in \psi^{-1}(\R^{d-d_0}\times W^{\prime})$. This in turn shows that $\Aa (\varphi(\yp),\yp) \in U$. Hence, $\Aa(\varphi(\yp),\yp) \in \Aa(G_{\varphi})\cap U = \Z \cap U$ and thus $g((\varphi(\yp),\yp)) = f(\Aa(\varphi(\yp),\yp)) = 0$.
Finally, using this in \cref{eq:intermediate}, recalling that $g = f\circ \Aa$, and defining $f_i: u \in U \mapsto g_i(\Aa^{-1}u)$, we have 
\begin{equation}
    \label{eq:finale}
    \forall u \in U,~ f(u) = \sum_{i=1}^{d-d_0}{f_i^2(u)}.
\end{equation}
We see that $f_i$ is of class $C^{p-2}$ since $g_i$ was of class $C^{p-2}$ and $\Aa^{-1}$ is an isometry; this concludes the proof of the theorem.
\end{proof}

\section{Discussion and possible extensions}\label{sec:discussion}\label{sec:extensions}\label{sec:conclusion}

In this work, we have provided second order sufficient conditions in order for a non-negative $C^p$ function to be written as a sum of squares of $C^{p-2}$ functions. We hope this will help provide a theoretical basis to algorithms which use functional sum of squares methods such as \cite{rudi2020finding,vacher21a,rudi2021psd}, which rely on the smoothness of such decompositions. As these conditions are sufficient and not necessary, one main problems is understanding this gap. This seems a highly difficult, and while very interesting, we present three other more reachable subjects for future work. 

The first is to have an explicit bound for the number of squares needed in the sum of squares decomposition in the compact case. We believe that using finer tools from differentiable geometry, we should be able to obtain a bound depending on meaningful quantities, and upper bounded by a constant $n_d$ depending only on the dimension $d$ of the manifold on which the function is defined.

The second is, as in the polynomial case, to handle functions $f$ which are non-negative on a constrained set defined by inequalities $f_i \geq 0$. More precisely, we would like to show second order sufficient conditions to write 
$f = g + \sum_{i}{g_i f_i}$ where the $g,g_i$ are sum of squares of regular functions when $f$ and the $f_i$ are regular. This would open up the field of constrained optimization for methods developed for functions, such as kernel sum of squares \cite{rudi2020finding}. In the polynomial case, such conditions are given by so-called \textit{Positivstellens\"atzen} (\cite{Putinar1993,Stengle1974,schmugden1991}), but usually assume the polynomial is positive. Second order conditions have been developed more specifically in \cite{marshall2006} to deal with non-negative polynomials with zeros. 

The third is to handle functions with conic outputs which are more general than the non-negativity one. For example, in order to represent functions which have values in a cone defined by linear inequalities   \cite{marteau2020non} or in the PSD cone  \cite{muzellec2022learning}. 
 
\bibliographystyle{siamplain}
\bibliography{references}

\appendix

\section{Around partitions of unity and gluing functions}\label{app:standard_manifolds}

In this section, we detail a few topological properties of manifolds, in order to \textit{a)} decompose a manifold or a sub-manifold in connected components and \textit{b)} use partitions of unity as a tool to glue functions together.
These specific properties are needed for \cref{sec:gluing}. 
For basics on topological spaces (what is a topology, the notion of continuity, of homeomorphism), we refer to Chapter 1 of \cite{topologie}. Main references for manifold can be found in \cite{lafontaine2015introduction,spivak,paulingeodiff}. Recall from \cref{sec:general_manifold} the definition of a manifold $M$ equipped with its atlas ${\cal A}$ of class $C^k$, and of a chart on $M$. A chart $\phi$ is said to be of class $C^{k^{\prime}}$ for $k^{\prime} \leq k$ if it compatible with the atlas up to $k^{\prime}$ smoothness, i.e. if the transitions maps $\phi \circ \phi_i^{-1}$ and $\phi_i \circ \phi^{-1}$ are all $C^{k^{\prime}}$. 
A priori, the atlas of a manifold of class $C^k$ is not unique in the sense that more than one atlas generate the same structure. To make it so, and to be able to say \textbf{the atlas} of $M$ of class $C^k$, we consider the \textbf{maximal atlas} on $M$, i.e. the collection of all charts of class $C^k$ on $M$.

\subsection{Paracompactness and partitions of unity}\label{app:paracompactness} The main point of asking a (differential) manifold to be second countable and Hausdorff, (and not just to be locally homeomorphic to $\R^d$), is for the manifold to be paracompact, and and hence to be equipped with partitions of unity. In this section, we introduce the main definitions and results on this topic.

Recall that a family of subsets $(U_{\alpha})$ of a space $X$ is said to be a covering of $X$ if $\bigcup_{\alpha}{U_{\alpha}} = X$. It is said to be locally finite if for any $x \in X$, there exists an open neighborhood $U$ of $x$ which intersects only a finite number of the $U_{\alpha}$. A family $(V_{\beta})$ is said to be a refinement of $(U_{\alpha})$ if for all $\beta$, there exists an $\alpha$ such that $V_{\beta} \subset U_{\alpha}$. 

A topological space $X$ is said to be \textbf{paracompact} if for any open covering $(U_{\alpha})$ of $X$, there exists an open refinement $(V_{\beta})$ of $(U_{\alpha})$ such that $(V_{\beta})$ is locally finite, and is an open covering of $X$.
The following lemma is proved in the first part of proposition 2.3 of \cite{paulingeodiff} or can be found in Theorem 2.13 of \cite{spivak}.
\begin{lemma}
A manifold is paracompact. 
\end{lemma}

Note that in \cite{spivak}, a manifold is defined to be a metric space locally like $\R^d$. In proposition 2.2 of \cite{paulingeodiff}, it is shown that being metric and second countable is equivalent to the countable Hausdorff condition (under the condition of being locally homeomorphic to $\R^d$). Spivak's condition in \cite{spivak} is however a bit more general; in fact, it allows a manifold $M$ to be a union of a possible non-countable connected component (as theorem 2 of \cite{spivak} shows that any connected component of a metric space locally homeomorphic to $\R^d$ is actually second countable).

Paracompactness is an important property as it yields the existence of partitions of unity. The following lemma is standard (a proof can be found in \cite{paulingeodiff}, proposition 2.3). The result is of course also true for $k = 0$, but is more technical to prove.

\begin{lemma}[Standard gluing lemma, \cite{paulingeodiff}]\label{lm:standard_gluing}
Let $(U_i)_{i \in I}$ be an open covering of a manifold $M$ of class $C^k$ (i.e. $\bigcup_{i \in I}{U_i} = M$). There exists a family of functions $\chi_i: M \rightarrow [0,1]$ of class $C^k$ such that $\supp(\chi_i) \subset U_i$ for all $i \in I$ and with locally finite support satisfying: 
\[\sum_{i \in I}{\chi_i} = 1.\]
\end{lemma}

We now prove the two technical results need in \cref{sec:gluing}.

\begin{proof}[Proof of \cref{lm:extension_lemma}]
The proof of this lemma is immediate. Indeed, by multiplication, we already know that $\chi g$ is well defined and $C^q$ on $U$. Moreover, for any point $x$ in $V = M \setminus{\supp(\chi)}$, which is an open set, $(\chi g)(x) = 0$ (by definition if $x \in M\setminus{U}$ and since $\chi(x) = 0$ if $x \in U$) and hence is $C^q$ on $V$. Since $V \cup U = M$ as $\supp(\chi)\subset U$, the property holds. Moreover, since $\chi g =0$ on $V$, $\supp(\chi g) \subset \supp(\chi)\subset U$.
\end{proof}

\begin{proof}[Proof of \cref{lm:glue_square}]
The proof of this result is a consequence of \cref{lm:standard_gluing}. Indeed, this result shows that there exists a family of function $\widetilde{\chi}_i: M \rightarrow [0,1]$ of class $C^k$ such that \textit{i)} for all $i \in I$, $\supp(\widetilde{\chi}_i)\subset U_i$, \textit{ii)} the support of $(\widetilde{\chi}_i)$ is locally finite and \textit{iii)} $\sum_{i}{\widetilde{\chi}_i} = 1$. 

Define $\phi = \sum_{i}{\widetilde{\chi}_i^2}$. Since $\sum_{i}{\widetilde{\chi}_i} = 1$, and $\widetilde{\chi} \geq 0$, necessarily, $\phi > 0$. Hence $\sqrt{\phi}$ is of class $C^k$, and hence $\chi_i:= \widetilde{\chi}_i / \sqrt{\phi}$ is of class $C^k$, and satisfies all the desired properties.
\end{proof}

\subsection{Connected components}

Connectedness is a key topological notion for manifolds, and allows to decompose manifold into separate blocks. Recall that two points $x,x^{\prime}$ of a topological set $X$ are connected if there exists no two open sets $U,V$ such that $X = U\cup V$, $x \in U$ and $x^{\prime} \in V$. Since being connected is an equivalence relation, we can partition $X$ in classes with respect to that relation, which are called "connected components". Connected components are both open and closed\footnote{For more details on connected components, see \cite{topologie}}. On a connected component of a manifold, the dimension $d$ of the charts $\phi : U \rightarrow \R^d$ is the same, and is called the dimension of that connected component (for more details, see any of the references on manifolds). 
Note that as a manifold $M$ is assumed to be second-countable, it has at most a countable number of connected components. Recall that a sub-manifold is defined in the main text as follows (such a definition can be found in section 2.4.2 of \cite{paulingeodiff}).

\begin{definition}\label{df:submanifold_manifold}
Let $M$ be a manifold of class $C^{k^{\prime}}$, $k \leq k^{\prime}$. $N$ is a sub-manifold of $M$ of class $C^k$ if for any $x \in N$, and any chart $\phi : U \rightarrow \R^d$ defined around $x$ and of class $C^k$, $\phi(U\cap N)$ is a sub-manifold (in the sense of $\R^d$, see \cref{secsec:submanifold}) of $\R^d$ around $\phi(x)$. It is equivalent to ask the existence of one such chart per point $x$.
\end{definition}

Let $N$ be a sub-manifold of class $C^k$ of a manifold $M$ of class $C^{k^{\prime}}$. Then it is naturally a manifold of class $C^k$ in its own right. Indeed, consider that \textit{i)} $N$ is equipped with the topology of $M$, i.e. $V$ is open in $N$ iif $V = U\cap N$ for some open set of $M$, and \textit{ii)} the atlas of $N$ is (the completion of) the set of restrictions of charts $\phi|_{U\cap N}$ where $\phi : U\rightarrow \R^d$ is a $C^k$ chart on $M$ such that $\phi(U \cap N) \subset \R^{d^{\prime}} \times \{0\}$, where $d^{\prime}$ is the dimension of $N$ at $x \in U$ (we identify $\R^{d^{\prime}} \times \{0\} \approx \R^{d^{\prime}}$). Note that the second-countable Hausdorff condition directly follows from that of $M$. Moreover, the $C^{k}$ compatibility of the charts is evident. From now on, when considering a sub-manifold $N \subset M$ as a manifold, it will be with this structure.
The reason for the introduction of all these concepts is to obtain the following lemma, which while it seems natural, we have not found as such in the literature.

\begin{lemma}\label{lm:open_sets}
Let $N$ be a sub-manifold of a manifold $M$. Let $(N_i)_{i \in I}$ be the connected components of $N$. There exists a collection of disjoint open sets $(U_i)_{i \in I}$ of $M$ such that each $N_i \subset U_i$.
\end{lemma}

\begin{proof} This proof relies mainly on paracompactness.

\noindent \textit{Step 1.} For all $x \in N$, there exists $U_x$ an open set in $M$ such that $\overline{U_x} \cap N$ is included in the unique connected component of $x$ in $N$. Indeed, by \cref{df:submanifold_manifold}, there exists a chart $\phi : U \rightarrow N$ where $U$ is an open neighborhood of $x$. But since $\phi(U\cap N)$ is a sub-manifold of $\R^d$ of class $C^k$ around $\phi(x)$, by Theorem 2.5 of \cite{paulingeodiff}, there exists a $C^{k}$ diffeomorphism $\psi : (\phi(x),V) \rightarrow (0,W)$ where $V$ such that $\psi(\phi(U\cap N)\cap V) = W \cap(\R^{d^{\prime}} \times \{0\})$ for some $d^{\prime}$. Taking $\widetilde{\phi} =\psi \circ \phi$ on $\widetilde{U} = \phi^{-1}(V)\cap U$, we have a chart of class $C^k$ around $x$ such that $\widetilde{\phi} : \widetilde{U} \rightarrow W\subset \R^d$ such that $\widetilde{\phi}(\widetilde{U} \cap N) = W \cap (\R^{d^{\prime}} \times \{0\})$. Now let $r$ be a radius such that the closed ball $\overline{B}(0,r) \subset W$. Set $U_x = \widetilde{\phi}^{-1}(B(0,r))$, which is an open neighborhood of $x$ included in $\widetilde{U}$. Note that $\overline{U_x} \subset  \widetilde{\phi}^{-1}(\overline{B}(0,r))\subset \widetilde{U}$ since $\widetilde{\phi}$ is continuous. Since $\overline{U}_x \cap N =  \widetilde{\phi}^{-1}(\overline{B}(0,r) \cap (\R^{d^{\prime}}\times \{0\}))$ which is connected, we have that $\overline{U}_x \cap N$ is connected and hence is included in the unique connected component of $x$.

\noindent \textit{Step 2.} Consider the collection of open sets $U_x$. By paracompactness of $U := \bigcup_{x \in N}U_x$ (it is a manifold), we can find an open cover $(U_{\alpha})$ of $U$ which is locally finite, and which still satisfies the condition that for all $\alpha$, $\overline{U}_{\alpha} \cap N$ is included in at most one connected component of $N$. Let $(N_i)_{ i \in I}$ denote the connected components of $N$. For $i \in I$, let $(V_{i,\alpha})_{\alpha \in A_i}$ denote the collection of open sets $U_{\alpha}$ such that $\overline{U}_{\alpha} \cap N \subset N_i$ and $\overline{U}_{\alpha} \cap N\neq \emptyset$. These collections satisfy \textit{a)} the $(V_{i,\alpha})_{\alpha \in A_i}$ cover $N_i$; \textit{b)} the collection $(V_{i,\alpha})_{i \in I,~\alpha \in A_i}$ has locally finite support; and \textit{c)} $\overline{V_{i,\alpha}}\cap N_i \subset N_i$ for all $i \in I,~\alpha \in A_i$.

\noindent \textit{Step 3.} For all $i \in I$, define $F_i =\bigcup_{j \in I\setminus{\{i\}},\beta \in A_j}{\overline{V}_{j,\beta}}$ and for all $\alpha \in A_i$, consider the set $W_{i,\alpha} = V_{i,\alpha} \setminus{F_i}$. $W_{i,\alpha}$ is open, and $W_{i,\alpha} \cap N = V_{i,\alpha}\cap N$. Indeed, let $x \in W_{i,\alpha}$. Since the $(V_{j,\beta})$ are locally finite, there exists $V_x \subset V_{i,\alpha}$ such that $V_x$ intersects a finite number of the $V_{j,\beta}$ and hence of the $\overline{V}_{j,\beta}$. Hence, $V_x \setminus{F_i}$ is still open. Hence $W_{i,\alpha}$ is open. The second condition comes from the fact that $\overline{V}_{i,\alpha} \cap N \subset N_i$, and that the connected components are disjoint
 Finally, taking $W_i = \bigcup_{\alpha \in A_i}{W_{i,\alpha}}$, the $W_i$ satisfy all the desired properties (they are disjoint thanks to the previous point and cover $N_i$ since the $V_{i,\alpha}$ covered $N_i$).
\end{proof}

\section{Morse lemma} 
In order for this article to be self contained, we restate the following classical lemmas from differential geometry and topology. 
Recall that a $C^k$-diffeomorphism is a map $\phi : U\subset \R^d \rightarrow V\subset \R^{d^{\prime}}$ which is of class $C^k$ and whose inverse is of class $C^k$ (in that case, necessarily, $d=d^{\prime}$). The following results are classical.
\begin{theorem}[Theorem 1.13 of \cite{lafontaine2015introduction}]\label{thm:loc_inversion}
Let $f: (x_0,\R^d) \rightarrow \R^d$ be a function of class $C^k$ ($k\geq 1$) defined around $x_0$ and such that $d f(x_0)$ is invertible. Then there exists a neighborhood $U$ of $x_0$ such that $f(U)$ is open and $f: U \rightarrow f(U)$ is a $C^k$ diffeomorphism.
\end{theorem}
\begin{theorem}[Theorem 1.18 of \cite{lafontaine2015introduction}]\label{thm:submersion}
Let $f: (x_0,\R^{d_1}) \rightarrow (y_0,\R^{d_2})$ be a function of class $C^k$ ($k\geq 1$) defined around $x_0$ s.t. $d f(x_0)$ is surjective and $f(x_0) = y_0$. Then there exists an open neighborhood $U$ of $x_0$ in $\R^{d_1}$, $V$ of $y_0$ in $\R^{d_2}$ as well as a function $g : V \rightarrow U$ of class $C^k$ such that $g(y_0) = x_0$ and $f\circ g = Id_{\R^{d_2}}$
\end{theorem}
We restate and reprove Lemma C.6.1 from \cite{Hormander2007}, which is a generalization of the so-called Morse Lemma (see lemma 2.2 of \cite{milnor1963}), and which is the basis of Morse Theory.
We will consider a function of two variables $f(x,y)$ defined on $\R^{d_1}\times \R^{d_2}$. We will denote with $\nabla_x f(x,y)$ its gradient with respect to the first variable taken at point $(x,y)$; it is an element of $\R^{d_1}$. Similarly, we will use the notation $\nabla^2_{xx}f(x,y) \in \R^{d_1\times d_1}$ to denote the Hessian matrix taken with respect to the first coordinate at point $(x,y)$. It is symmetric.

\begin{lemma}[Lemma C.6.1 from \cite{Hormander2007}] \label{lm:morse_hormander}
Let $d_1,d_2 \in \N, p \in \N$ with $p \geq 2$. Let $f: (x,y) \in U_0 \subset \R^{d_1} \times \R^{d_2} \mapsto f(x,y) \in \R$ be $C^p$ function defined on a neighborhood $U_0$ of $(0,0)$. Assume that $\nabla_{x}f(0,0) = 0$ and that $H:= \nabla^2_{xx}f (0,0)$ is non-singular.  

There exists an open convex neighborhood $V$ of $0$ in $\R^{d_1}$ and an open convex neighborhood $W$ of $0$ in $\R^{d_2}$ such that $V \times W \subset U_0$, a map $\varphi \in C^{p-1}(W,V)$ and a map $z \in C^{p-2}(V \times W,\R^{d_1})$ such that for any $(x,y) \in V \times W$  $\nabla_x f(x,y) = 0$ if, and only if $x = \varphi(y)$, and 
\begin{equation}
\label{eq:morse_lemma_decomposition}
\forall (x,y) \in V\times W,~f(x,y) = f(\varphi(y),y) + \frac{1}{2}z(x,y)^\top H z(x,y).
\end{equation}
\end{lemma} 

To simplify the proof, we first show an intermediate result which gives $\varphi$. 

\begin{lemma}\label{lm:technical_hoho}
Under the assumptions of \cref{lm:morse_hormander}, there exists two open convex neighborhoods of zero $V_0 \subset \R^{d_1},~W_0\subset \R^{d_2}$ and $\varphi: W_0 \rightarrow V_0$ of class $C^{p-1}$ such that \textit{a)}  $V_0\times W_0 \subset U_0$ and \textit{b)} $\forall (x,y) \in V_0\times W_0,~\nabla_{x}f(x,y) = 0 \Leftrightarrow x = \varphi(y)$.
\end{lemma}
\begin{proof}
Consider the map $\psi: (x,y) \in U_0 \subset \R^{d_1} \times \R^{d_2} \mapsto (\nabla_{x}f(x,y),y)$. Its jacobian at $(0,0)$ is of the form  $\begin{pmatrix}H & \star\\ 0& \Id_{d_2}\end{pmatrix}$. 
Since $H$ is non-singular, this matrix is non-singular. Applying the local inversion lemma \cref{thm:loc_inversion}, there exists an open neighborhood $U_1 \subset U_0$ such that $\psi$ is a $C^{p-1}$ diffeomorphism from $U_1$ to  $\psi(U_1)$.

Let $\widetilde{V}_0 \subset \R^{d_1}$, $\widetilde{W}_0 \subset \R^{d_2}$ be open convex neighborhoods of $0$ such that $\widetilde{V}_0 \times \widetilde{W}_0 \subset U_1 \cap \psi(U_1)$. Define $\varphi: w \in \widetilde{W}_0 \mapsto \pi_{1}(\psi^{-1}(0,w)) \in \R^{d_1}$. Defining $V_0 = \tilde{V}_0$ and $W_0 \subset \R^{d_2}$ to be an open convex neighborhood of $0$ included in $\varphi^{-1}(V_0) \cap \widetilde{W}_0$, we have $\varphi(W_0) \subset V_0$ and $V_0 \times W_0 \subset U_1 \subset U_0$.

 Moreover, for any $(x,y) \in V_0 \times W_0 \subset U_1 \cap \psi(U_1)$, $\nabla_1 f(x,y) = 0$ iif $\psi(x,y) = (0,y) \in \psi(U_1)$, iif $(x,y) = \psi^{-1}(0,y) = (\varphi(y),y)$, iif $x = \varphi(y)$.
\end{proof}
We can now prove our main result.
\begin{proof}[Proof of \cref{lm:morse_hormander}]
 Fix $V_0,W_0$ satisfying the properties of \cref{lm:technical_hoho}. Let $(x,y)\in V_0 \times W_0$. For $t \in [0,1]$, define $x_t = \varphi(y) + t(x-\varphi(y))$. By convexity of $V_0$, $(x_t,y) \in V_0 \times W_0 \subset U_1 \subset U_0$ for all $t \in [0,1]$. Thus, the map $g: t \in [0,1] \mapsto f(x_t,y)$ is well defined, and we can apply the Taylor formula $g(1) = g(0) + g^{\prime}(0) + \int_{0}^1{(1-t)g^{\prime\prime}(t)dt}$ and the fact that $g^{\prime}(0) = \nabla_x f(\varphi(y),y)\cdot (x-\varphi(y)) = 0$ to obtain
\[f(x,y) = f(\varphi(y),y) + (x - \varphi(y))^\top \left(\int_{0}^1{(1-t)\nabla^2_{xx}f(x_t,y)dt}\right)(x-\varphi(y))\]
Defining $B: V_0 \times W_0 \rightarrow S(\R^{d_1})$, such that $B(x,y):= 2\int_{0}^1{(1-t)\nabla^2_{xx}f(x_t,y)dt}$, the previous equation  can simply be written $f(x,y) = f(\varphi(y),y) + \tfrac{1}{2} (x-\varphi(y))^\top B(x,y) (x-\varphi(y))$. Note that
$B \in C^{p-2}(V_1 \times W_1,S(\R^{d_1}))$ and $B(0,0) = H$.
Now define $G: R \in \R^{d_1 \times d_1} \mapsto R^\top H R \in S(\R^{d_1})$ which is $C^{\infty}$ and whose differential in $\Id_{\R^{d_1}}$ is surjective (see \cite{Hormander2007}). \cref{thm:submersion} shows there exists an neighborhood ${\cal O}$ of $H$ in $S(\R^{d_1})$ and a $C^{\infty}$ function $F: {\cal O} \rightarrow  \R^{d_1 \times d_1}$ such that $(G \circ F)(B) = B$ for all $B \in {\cal O}$. 
Let $V \subset \R^{d_1},\widetilde{W} \subset \R^{d_2}$ be two open convex neighborhoods of $0$ such that $V \times \widetilde{W} \subset B^{-1}(\cal O)$. Let $W$ be an open convex neighborhood of $0$ such that $W \subset \widetilde{W} \cap \varphi^{-1}(V)$ and define $z(x,y) = (F\circ B)(x,y)(x-\varphi(y))$. $z$ satisfies \cref{eq:morse_lemma_decomposition}. 
\end{proof}

\end{document}

%% file: siopt_shared.tex
\usepackage{lipsum}
\usepackage{amsfonts}
\usepackage{graphicx}
\usepackage{epstopdf}
\usepackage{algorithmic}

\usepackage{booktabs}
\usepackage{color}
\usepackage{latexsym} 
\usepackage{amsmath}               
\usepackage{amssymb}               
\usepackage{tcolorbox}
\usepackage{times}
\usepackage[utf8]{inputenc}
\usepackage[T1]{fontenc}
\usepackage{enumitem}
\usepackage{cite}
\usepackage{nicefrac}       
\usepackage{microtype}      
\usepackage[american]{babel}
\usepackage{thmtools}
\usepackage{booktabs}
\usepackage[tikz]{bclogo}
\usepackage{multicol,tabularx,capt-of}
\usepackage{hhline}
\usepackage{multirow}

\usepackage{microtype}

\usepackage{stmaryrd}
\usepackage{hyperref}
\usepackage{xcolor}
\definecolor{dgreen}{rgb}{0.00,0.49,0.00}
\definecolor{dblue}{rgb}{0,0.08,0.75}
\hypersetup{
    colorlinks = true,
    citecolor=dgreen,
    urlcolor=dblue,
    linkcolor=dblue
}

\usepackage[parfill]{parskip}

\usepackage{subfigure} 

 \usepackage{cite}
\usepackage{newtxtext}      

\newsiamremark{remark}{Remark}
\newsiamremark{example}{Example}
\crefname{example}{Example}{Examples}
\newsiamremark{hypothesis}{Hypothesis}
\crefname{hypothesis}{Hypothesis}{Hypotheses}
\newsiamthm{claim}{Claim}
\crefname{theorem}{Theorem}{Theorems}

\headers{Decomposing smooth functions as sums of squares}{Ulysse Marteau-Ferey, Francis Bach, Alessandro Rudi}

\title{Second order conditions to decompose smooth functions\\ as sums of squares \thanks{This work was funded in part by the French government under management of Agence Nationale de la Recherche as part of the “Investissements d’avenir” program, reference ANR-19-P3IA-0001(PRAIRIE 3IA Institute). We also acknowledge support from the European Research Council (grants SEQUOIA 724063 and REAL 947908), and support by grants from Région Ile-de-France.}}

\author{Ulysse Marteau-Ferey\thanks{Inria/PSL Research University
  (\email{ulysse.marteau-ferey@inria.fr},\email{francis.bach@inria.fr},\email{alessandro.rudi@inria.fr}).}
\and Francis Bach \footnotemark[2]
\and Alessandro Rudi \footnotemark[2]
}

\usepackage{amsopn}

\newcommand{\R}{\mathbb{R}}

\newcommand{\N}{\mathbb{N}}
\newcommand{\Z}{{\cal Z}}

\newcommand{\cc}{{\cal C}}

\newcommand{\Id}{\boldsymbol{I}}
\newcommand{\xpar}{x_{\shortparallel}}
\newcommand{\xperp}{x_{\perp}}

\newcommand{\psd}{\mathbb{S}_+}

\newcommand{\xp}{x^{\prime}}
\newcommand{\yp}{y^{\prime}}

\newcommand{\Aa}{{\cal A}}

\newcommand{\Wt}{\widetilde{W}}

\DeclareMathOperator{\Hom}{Hom}

\DeclareMathOperator{\im}{im}
\DeclareMathOperator{\supp}{supp}
